\documentclass[10pt]{amsart}
\usepackage{amsfonts} 
\usepackage{amscd}
\usepackage{amssymb}
\usepackage{amsrefs}
\input xy
\xyoption{all}

\newcommand{\eg}{{\it{e.g.}}}

\newcommand{\ring}[1]{\mathbb{#1}}
\newcommand{\Cgordon}{\ring{C}} \newcommand{\Qgordon}{\ring{Q}}
\newcommand{\Zgordon}{\ring{Z}} \newcommand{\R}{\ring{R}}
\newcommand{\A}{\ring{A}}
\newcommand{\F}{\ring{F}} 
\newcommand{\N}{\ring{N}}

\newcommand{\lef}{\ring{L}}

\newcommand{\be}{\begin{equation}}
\newcommand{\ee}{\end{equation}}

\newcommand{\ta}[1]{{#1}[[t]]}
\newcommand{\arcs}{\mathfrak{L}}
\newcommand{\mor}{\mathrm{Mor}}

\newcommand{\mot}{{\text{Mot}}}

\newcommand{\spec}{\mathnormal{\mathrm{Spec\,}}}
\def\1{{\mu\mkern-6mu\mu}}
\newcommand{\lie}{{\bf\mathfrak g}}

\newcommand{\Ggordon}{{\mathbb G}}

\newcommand{\ord}{{\text{ord}}}

\newcommand{\rigordon}{{\mathcal O}}

\newcommand{\pff}{{\text{PFF}}}
\newcommand{\dt}{\,{\text d}t}
\newcommand{\ie}{{\it i.e.}}
\newcommand{\jac}{{\text{Jac}}}
\newcommand{\var}{{\text{Var}}}
\newcommand{\epsgordon}{{\varepsilon}}
\newcommand{\calf}{{\mathcal C}}
\newcommand{\id}{{\operatorname d}}

\newcommand{\rf}{\F}
\newcommand{\valgordon}{{\ord}}

\newcommand{\ac}{\overline{\text{ac}}}
\newcommand{\mx}{{\mathcal X}}

\newcommand{\dergordon}{{\text {d}}}

\newcommand{\de}{{\text{Def}}}
\newcommand{\rde}{{\text{RDef}}}

\theoremstyle{plain}
\newtheorem{thm}{Theorem}
\newtheorem{theorem}[thm]{Theorem}
\newtheorem{proposition}[thm]{Proposition}

\theoremstyle{definition}
\newtheorem{remark}[thm]{Remark}
\newtheorem{definition}[thm]{Definition}
\newtheorem{example}[thm]{Example}
\newtheorem{exercise}[thm]{Exercise}
\author[]{Julia Gordon and  Yoav Yaffe}

\title[]
{An overview of arithmetic motivic integration}

\begin{document} 
%\chapter{An overview of arithmetic motivic integration}{Julia Gordon}  

\maketitle

\section{Introduction}

The aim of these notes is to provide an elementary
introduction to some aspects of the theory of 
arithmetic motivic integration, as well as a brief guide to the extensive 
literature on the subject. 

The idea of motivic integration was introduced by M. Kontsevich in 1995.
It was quickly developed by J. Denef and F. Loeser in a series of papers 
\cite{DL}, \cite{DL.McKay}, \cite{DL.Igusa}, and {by} others.
This theory, whose applications are mostly in algebraic geometry over 
algebraically closed fields,  now is often referred to as 
``geometric motivic integration'', 
to distinguish it from the so-called arithmetic motivic integration that 
specifies to integration over $p$-adic fields.

The theory of arithmetic motivic integration first appeared in the 1999 paper by J. Denef 
and F. Loeser  \cite{DL.arithm}. 
The articles \cite{Tom.intro} and \cite{DL.congr} together provide  an
excellent exposition of this work. 
In 2004, R. Cluckers and F. Loeser developed a different and very effective
approach to motivic integration (both geometric and arithmetic) \cite{CL}.
Even though there is an expository version 
\cite{CL.expo}, this theory seems to be not yet well-known. 
This note is intended in part to be a companion with examples 
to \cite{CL}. The aim is not just to describe 
what motivic integration achieves, but to give some clues as to how it works.
We have stayed very close to the work of Cluckers and Loeser in the 
main part of this exposition. In fact, much of these notes is a direct 
quotation, most frequently from
the articles \cite{CL}, \cite{CL.expo},  and also 
\cite{DL.congr}, and \cite{DL.arithm}. 
Even though we try to give precise references all the time, some quotes from these sources might not always be acknowledged since they are so ubiquitous. 
Some ideas, especially in the appendices, are clearly borrowed from 
\cite{Tom.intro}. 
The secondary goal was to collect  references to many sources on motivic integration, and to provide some information on the relationship and logical interconnections between them. This is done in Appendix 1 (Section~\ref{geo}). 

Our ultimate hope is that the reader would be able to start using motivic integration instead of $p$-adic integration, if there is any advantage in doing 
integration independently of $p$ at the cost of losing a finite number of primes.

{\bf Acknowledgment.}
The first author thanks T.C. Hales for introducing her to the subject; 
Jonathan Korman -- for many hours of discussions, and 
Raf Cluckers -- for explaining his work on several occasions.
We have learned a lot of what appears in these notes at the 
joint University of Toronto-McMaster University seminar on motivic 
integration in 2004-2005, and thank  Elliot Lawes, Jonathan Korman and 
Alfred Dolich for their lectures. 
{The contributions of
 the second author are limited to sections 1-5.}
Finally, the first author thanks the organizers and participants 
of the mini-courses  on motivic integration at the University of Utah and at 
the Fields Institute Workshop at the University of Ottawa, where most of 
this material was presented, and the editors of this volume for multiple 
suggestions and corrections.

%\tableofcontents

\section{$p$-adic integration}
Arithmetic motivic integration re-interprets the classical 
measure on $p$-adic fields, and $p$-adic manifolds, in a geometric way.
The main benefit of such an interpretation is that it allows one to isolate the dependence on $p$, so that 
one can perform integration in a field-independent way, and then ``plug in'' $p$ at the very end.
Even though this is not the only achievement of the theory, it will be our main focus in these notes.
Hence, we begin with a brief review of the properties of {the} field of 
$p$-adic numbers, and integration on $p$-adic manifolds. 

\subsection{The $p$-adic numbers}
Let $p$ be a fixed prime.
Throughout these notes our main example {of a local field} will be the field $\Qgordon_p$ of $p$-adic numbers, which is the completion of $\Qgordon$ with respect to the $p$-adic metric.
\subsubsection{Analytic definition of the field $\Qgordon_p$}
Every non-zero rational number $x\in\Qgordon$ can be written in  the form 
$x=\frac{a}{b}p^n$, where $n\in\Zgordon$, and $a,b$ are integers relatively 
prime to $p$. The power $n$ is called the {\bf valuation} of $x$ and 
denoted $\ord(x)$. Using the valuation map, we can define a
norm on $\Qgordon$: $|x|_p=p^{-\ord(x)}$ if $x$ is non-zero and 
{$|0|_p = 0$.} This norm induces a metric on $\Qgordon$, which satisfies a stronger triangle inequality than the standard metric:
$$|x+y|_p\le \max\{|x|_p, |y|_p\}.$$
This property of the metric is referred to as the ultrametric property.

The set $\Qgordon_p$, as a metric space, is the completion of $\Qgordon$ with respect to this metric. The operations of addition and multiplication extend by continuity
from $\Qgordon$ to $\Qgordon_p$ and make it a field.
The set $\{x\in \Qgordon_p\mid \ord(x)\ge 0\}$ is denoted $\Zgordon_p$ and called the ring of $p$-adic integers.

\subsubsection{Algebraic definition of the field $\Qgordon_p$}
There is a way to define $\Qgordon_p$ without invoking analysis. 
Consider the {rings} $\Zgordon/p^n\Zgordon$. They form a projective system with natural maps
$$
\begin{aligned}
\Zgordon/p^{n+1}\Zgordon&\to\Zgordon/p^n\Zgordon \\
m&\mapsto m \pmod{p^n}.
\end{aligned}
$$
The projective limit is called $\Zgordon_p$, the ring of $p$-adic integers.
The field $\Qgordon_p$ is then defined to be its field of fractions.

\subsubsection{Basic facts about $\Qgordon_p$}
\begin{itemize}

\item The two definitions of $\Qgordon_p$ agree, and $\Qgordon_p$ is a field extension of
$\Qgordon$.

\item
Topology on $\Qgordon_p$: if we use the analytic definition, then $\Qgordon_p$ comes  equipped with a metric topology. It follows from the strong triangle inequality that $\Qgordon_p$ is totally disconnected in this topology. 
It is easy to prove  that the sets $p^n\Zgordon_p$, as $n$ ranges over $\Zgordon$, form a basis of neighbourhoods of $0$.
If one uses the algebraic definition of $\Qgordon_p$, then 
the topology for $\Qgordon$ is \emph{defined} by declaring that these sets form a basis of neighbourhoods of $0$, and
the basis of neighbourhoods at any other point is obtained by translating them. 

\item The set $\Zgordon_p\subset \Qgordon_p$ is open and compact in this topology.
It follows that each $p^n\Zgordon_p$ is also a compact set, 
{ which, in turn, implies that} $\Qgordon_p$ is 
locally compact.
Note that $\Zgordon_p$ (in the analytic definition) has the description 
$\Zgordon_p=\{x\in \Qgordon_p\mid |x|_p\le 1\}$, so it is the {closed} 
unit ball in our metric 
space (somewhat counter-intuitively).
 
\item As a set, $\Zgordon_p$ is in bijection with the set 
$$\left\{\sum_{i=0}^{\infty}a_i p^i\mid a_i=0,\dots, p-1\right\}.$$
Note, however, that addition in $\Zgordon_p$ \emph{does not} agree with 
coefficient-wise addition $\mod p$ of the power series 
(because ``$p$ has to carry over'').\footnote{Passing to the 
fields of fractions, we see that the field $\F_p((t))$ of formal Laurent series with 
coefficients in $\F_p$, and $\Qgordon_p$ are naturally in bijection, but 
not isomorphic;  
we will see that these fields have a lot in common nevertheless.}

\end{itemize}

In these notes, we work with discretely valued fields, \ie, 
fields equipped with a valuation map from the non-zero elements of the field to a group $\Gamma$ with a discrete topology; this valuation will always be denoted by $\ord$.
From now on, we will always assume that $\Gamma=\Zgordon$.

\begin{theorem}\label{thm:classif}
Any complete discretely valued field that is locally compact in the topology induced by the valuation is isomorphic either to a finite extension of $\Qgordon_p$ or to a field $\F_q((t))$ of formal Laurent series over a finite 
field.  
\end{theorem}

We refer to fields of this kind by the term ``local fields''; when we want to distinguish between finite extensions of $\Qgordon_p$ and the function fields $\F_q((t))$, we refer to them as ``characteristic zero fields'', and 
``equal characteristic fields'', respectively.

\begin{remark} Note that if a field $k((t))$ of formal 
Laurent series over a field $k$ is locally compact, then $k$ is finite.
\end{remark}

The above theorem and a discussion of related topics can be found, 
for example, in \cite[Appendix to Chapter 2]{Robert}.

\subsection{Hensel's Lemma}
{The theory of integration on local fields that is the focus of these notes} would have been impossible without the property of the non-archimedean local fields known as Hensel's Lemma.
The next example is classical; we include it as a reminder.

\begin{example} $\Zgordon_3$ does not contain $\sqrt{2}$, but $\Zgordon_7$ does. Indeed,
let us try to solve the equation $x^2=2$. If we write 
$x=\sum_{i=0}^{\infty} a_i 3^i$, {then} $x^2=a_0^2+3\cdot 2a_0a_1 +3^2(a_1^2+2a_0a_2)+\dots$, {hence}  $x^2\pmod{3}=a_0^2 \pmod{3}$.  Since $a_0^2$ cannot 
be congruent to 
$2 \pmod{3}$, there is no solution.

However, if we play the same game $\mod 7$, we {have solutions, for example} $a_0=3$. Next we need to find $a_1$ such that $(3+7a_1)^2\equiv 2 \pmod{49}$ (we find $a_1=1$), and so on.
{Clearly for every step $i\ge 1$ we can find a unique solution for $a_i$,} and this way we get a power series which converges 
(in $\Qgordon_7$) to a solution of the equation $x^2=2$. Since $\Qgordon_7$ is complete, it must 
contain the sum of the series, and this way $\sqrt{2}$ is in $\Qgordon_7$. Since the series has no negative powers, it is {in $\Zgordon_7$, but this actually follows from $|\sqrt{2}|_7 = \sqrt{|2|_7}=\sqrt{1}=1$.}
\end{example}

\begin{theorem}(Hensel's Lemma)
Let $K$ be a non-archimedean local field, and let $f\in K[x]$ be a monic polynomial such that all its coefficients have non-negative valuation. Then if $x\in K$ has the property that $\ord(f(x))>0$ and $\ord(f'(x))=0$, then there exists
$y\in K$ such that $f(y)=0$ and $\ord(y-x)>0$.
\end{theorem}

The root $y$ is constructed by using Newton's approximations (and taking $x$ as the first one). The completeness of the field  is used to establish 
convergence of the sequence of appoximations to a root $y$.

The meaning of Hensel's Lemma (\eg\ for $\Qgordon_p$) is that every solution 
of $f(x)=0 \pmod{p}$ can be lifted to an actual root of $f(x)$ in $\Zgordon_p$
(in particular, it justifies our example of $\sqrt{2}\in \Zgordon_7$).

Fields satisfying Hensel's Lemma are called Henselian.
The argument sketched above shows that all complete discretely valued 
fields are Henselian, in particular,  the fields of formal Laurent series
$K((t))$ where $K$ is an arbitrary field, are Henselian.
See \eg \cite{ribenboim} for a detailed discussion of the Henselian property.  

\subsection{Haar measure}\label{sub:haar}
{If $K$ is a locally compact non-archimedean field then the additive group 
of  $K$ (as a locally compact abelian group) has a  unique up to a constant multiple translation-invariant 
measure, called the Haar measure.} 
The 
$\sigma$-algebra of measurable sets is the usual Borel $\sigma$-algebra (generated by open sets in the topology on $K$ induced by the absolute value on $K$). This measure will be denoted  by $\mu$.

It is easy to check that $\mu$ satisfies a natural 
``Jacobian rule'': if $a\in K$, and $S\subset K$ is a measurable set, then  
$\mu(aS)=|a|\mu(S)$.

\begin{example} Though it is very simple, this example is the source of intuition behind much of our general theory: if we normalize the Haar measure on $\Qgordon_p$ so that $\mu(\Zgordon_p)=1$, then
$\mu(p^n\Zgordon_p)=p^{-n}$.
\end{example}

Using the product measure construction, we can get a translation-invariant 
measure on the affine space $\A^n(K)$ from the Haar measure on $K$. It is also unique up to a constant multiple.
Since it plays an important role in the construction of motivic measure, we recall Jacobian transformation rule for the $p$-adic measure, which is analogous 
to the transformation rule for Lebesgue measure on $\R^n$. 

\begin{theorem} Let $A$ be a measurable subset of $\A^n(\Qgordon_p)$, let $\phi$ 
be a 
$C^1$-map $\phi:\A^n\to \A^n$  injective and with nonzero Jacobian on $A$, 
and let 
$f:\A^n(\Qgordon_p)\to \R$ be an integrable function.
Then 
$$
\int_A f\dergordon \mu=\int_{\phi(A)}|\jac(\phi)|_pf\circ\phi^{-1}\dergordon \mu.
$$
\end{theorem}

\subsection{Canonical (Serre-Oesterl\'e) measure on $p$-adic 
manifolds}\label{sub:Weil}
Much of this section is quoted from \cite{batyrev}; see also \cite{Weil}.

Let $K$ be a local Henselian field with valuation $\ord$, uniformizer 
$\varpi$, the ring of integers ${\mathcal O}_K=\{x\in K\mid \ord(x)\ge 0\}$, and the residue field $\F_q$.
Let ${\mathcal X}$ be a smooth scheme over $\spec\rigordon_K$ of 
dimension $d$.\footnote{Alternatively, one can think of ${\mathcal X}$ as
a smooth variety over $K$, such that its reduction $\mod \varpi$ is a smooth 
variety  over $\F_q$. The smooth subschemes $U_i$ below then should be 
replaced with  Zariski open subsets.} 
Assume for now that there is a nowhere vanishing global differential form 
$\omega$ on 
${\mathcal X}$.
Since ${\mathcal X}$ is smooth, ${\mathcal X}(\rigordon_K)$ is a $p$-adic manifold, and we can use {the} differential form $\omega$
to define a measure on ${\mathcal X}(\rigordon_K)$, in the following way.

Let $x\in {\mathcal X}(\rigordon_K)$ be a point, and $t_1,\dots,
t_d$ be the local coordinates around this point.
They define a homeomorphism $\theta$
(in {the} $p$-adic analytic topology) from an open neighbourhood 
$U\subset {\mathcal X}(\rigordon_K)$ to an open set $\theta(U)\subset \A^d(\rigordon)$.
We write
$$
\omega=\theta^{\ast}(g(t_1,\dots, t_d)\dergordon t_1\wedge \dots \wedge\dergordon t_n),
$$
where $g(t_1,\dots, t_n)$ is a $p$-adic analytic function on $\theta(U)$ 
with no zeroes.
Then we can define a measure on $U$ by $\dergordon\mu_{\omega}=|g(t)|\dergordon t$, where
$|\dergordon t|$ stands for the measure on $\A^d$ associated with the volume form 
$\dergordon t_1\wedge\dots \wedge \dergordon t_d$ (\ie, the product measure defined in the previous subsection), normalized so that 
$$\int_{\A^d(\rigordon_K)}|\dergordon t|=1.$$
Two different nonvanishing differential forms on ${\mathcal X}$ have to differ by a $p$-adic unit; therefore, they yield the same measure on 
${\mathcal X}(\rigordon_K)$. 

More generally, even if there is no non-vanishing form, we can cover 
${\mathcal X}$ with finitely many smooth affine open subschemes $U_i$
such that on each one of
them there is a nonvanishing top degree form $\omega_i$. 
The form $\omega_i$ allows us to transport the measure from $\A^d(\rigordon_K)$
to $U_i(\rigordon_K)\subset {\mathcal X}(\rigordon_K)$. Note that each of the forms $\omega_i$
is defined uniquely up to an element 
$s_i\in\Gamma(U_i,\rigordon_{\mathcal X}^{\ast})$.
Therefore, the measure we
define
on $U_i(\rigordon_K)$ does not depend on the choice of $\omega_i$, since
by definition, $|s_i(x)|=1$ for   
$s_i\in\Gamma(U_i,\rigordon_{\mathcal X}^{\ast})$, $x\in U_i$.
These measures on $U_i(\rigordon_K)$ glue together (to check this, we need
to check that our measures on $U_i(\rigordon_K)$ and $U_j(\rigordon_K)$ agree on
the overlap $U_i(\rigordon_K)\cap U_j(\rigordon_K)=U_i\cap U_j(\rigordon_K)$. This
follows from the definition of a differential form and the fact that
the measure on $\A^d(\rigordon_K)$ satisfies the Jacobian transformation
rule.

\begin{definition} The measure defined as above on 
${\mathcal X}(\rigordon)$ is 
called {\bf the canonical $p$-adic measure}.
\end{definition} 

The above approach (due to A. Weil) 
to the definition of the measure has, through the 
work of Batyrev, inspired the initial approach to motivic integration.
{Let us also sketch Serre's definition of the canonical measure on 
subvarieties of the affine space, \cite{serregordon}, which was generalized 
by Oesterl\'e to $p$-adic analytic sets \cite{Oesterle}, and by W. Veys -- to subanalytic sets, 
\cite{Veys:measure},} 
and which is now used as the classical definition of the $p$-adic measure.
Let $Y$ be a $d$-dimensional smooth subvariety of $\A^n$.
Instead of using local coordinates, one can use the coordinates of the 
ambient affine space $x_1,\dots, x_n$. 
For each subset of indices $I=\{i_1,\dots, i_d\}$ with $i_1<i_2<\dots<i_d$,
let $\omega_{Y,I}$ be the differential form on $Y$ induced by 
$dx_{i_1}\wedge\dots\wedge dx_{i_d}$.
Let $\mu_{Y,I}$ be the measure on $Y$ associated with the form $\omega_{Y,I}$.
The canonical measure on $Y$ is defined by
$$\mu_Y=\sup_{I}\mu_{Y,I},
$$ 
where $I$ runs over all $d$-element subsets of $\{1,\dots, n\}$.

\subsection{Weil's theorem}

\begin{theorem}(A. Weil, \cite{Weil}) Let $K$ be a locally compact
  non-archimedean field with the ring of integers $\rigordon_K$ and residue
  field $\F_q$, and let ${\mathcal X}$ be a smooth scheme over $\spec
  \rigordon_K$ of dimension $d$. Then
$$\int_{{\mathcal X}(\rigordon_K)}\dergordon \mu=\frac{|\mathcal X(\F_q)|}{q^d}.
$$
\end{theorem}

\noindent{\bf Heuristic ``proof''.}
Consider the projection from ${\mathcal X}(\rigordon)$ to 
${\mathcal X}( \rf_q)$ that is defined by 
applying the  reduction $\mod \varpi$ map to the local coordinates. 
The smoothness of ${\mathcal X}$ implies that the fibres 
of this map all look exactly like the fibres of this projection for the affine 
space.
The definition of the measure on ${\mathcal X}(\rigordon)$ implies that 
the measure on ${\mathcal X}(\rigordon)$ is transported from the measure on 
the affine space, on each coordinate chart.  
Finally, in the case of the affine space of dimension $d$, 
the volume of each fibre of the projection is $q^{-d}$.
Hence, the total volume of ${\mathcal X}(\rigordon)$ equals the cardinality of the 
image of the projection times the volume of the fibre, which is,
$|{\mathcal X}(\rf_q)|q^{-d}$.

Motivic integration initially started (in a lecture by M. Kontsevich) 
as a generalization of these ideas behind $p$-adic integration.\footnote{Here we have (rather carelessly and briefly) considered only 
the $p$-adic manifolds that arise by taking $\Zgordon_p$-points of smooth varieties.
It should at least be mentioned that there is a theory of 
motivic integration on rigid analytic spaces \cite{Sebag}.}

\subsection{Motivation}
The goal of arithmetic motivic integration is to assign a geometric object 
to a 
$p$-adic measurable set, in such a way that the value of the measure can 
be recovered by counting the number of points of this geometric object over 
the residue field, in a way that is similar to Weil's Theorem, 
where the volume of ${\mathcal X}(\rigordon_K)$ is recovered by counting 
points on the closed fibre of ${\mathcal X}$, for a smooth projective 
scheme ${\mathcal X}$.

There are two approaches to the development of arithmetic motivic integration.
The first one, that appears in \cite{DL.arithm}, uses arc spaces and 
truncations. 
The missing link between the residue fields of characteristic zero and residue fields of finite characteristic is provided by considering pseudofinite fields. There is a beautiful exposition of this 
work \cite{Tom.intro}. We sketch the main steps in Appendix 1, for 
completeness of this overview.
  
The other approach uses cell decomposition instead of truncation. 
This theory was developed in \cite{CL}, and it gives more than just a 
measure -- it is a theory of integration complete with a large class of 
integrable functions. 
The main body of these notes is devoted to this theory. 
In order to describe it, we need to introduce some techniques from logic and 
some abstract formalism. This is done in the next section. For a while there will be no $p$-adic manifolds, and no measure. These familiar concepts will 
reappear in Section~\ref{sec:back}.

\section{Constructible motivic Functions}\label{sec:3lang}
To start with, we need to develop a way to talk about sets without mentioning their elements, similar to the way a variety is defined independently of its set of points over any given field. 
This is done by means of specifying a language of logic, and describing sets by
formulas in that language.

Given a formal language $L$, a we say that an object $M$ is a 
{\bf structure} for this language if formulas in $L$ can be interpreted 
as statements about the elements of $M$. More precisely, 
in this context one has an interpretation function  from the set of 
symbols of the  alphabet of $L$ to $M$, a map from the symbols for 
operations in $L$ to $M$-valued functions on the direct product of the corresponding number of copies of $M$, and a map that takes symbols for relations
in $L$ (such as ``='', for example) to subsets of $M^r$ with the 
corresponding $r$ (\eg, $r=2$ for binary relations such as ``='').  
We refer for example to \cite{manin.logic} for a discussion of this and related topics.

Let $L$ be a language, and let $M$ be a structure for $L$.
A set $A\subset M^n$ is called {\bf $L$-definable} if there exists a formula in the language $L$ such that $A$ is the set of points in $M^n$ satisfying this 
formula.

A function is called {\bf $L$-definable} if its graph is a definable set.   

\subsection{The language of rings}\label{sub:ldp}
The first-order language of rings is the language {built from the 
following set of symbols}:
\begin{itemize}
\item countably many symbols for variables $x_1, \dots, x_n,\dots$.
\item '$0$','$1$';
\item '$\times$', '$+$', '$=$', and parentheses '$($', '$)$';
\item The existential quantifier '$\exists$';
\item logical operations: conjunction '$\wedge$', negation '$\neg$', disjunction'$\vee$'.
\end{itemize}

Any syntactically correct formula built from these symbols is a formula in 
the first order language of rings.
 
Any ring is a structure for this language. 

Note that quantifier-free 
formulas in the language of rings define constructible sets
(recall that constructible sets, by definition, are the sets that belong to the smallest family ${\mathcal F}$ containing Zariski open sets and such that
a finite intersection of elements of ${\mathcal F}$ is in ${\mathcal F}$, and a complement of an element of ${\mathcal F}$ is in ${\mathcal F}$). 

\subsection{Presburger's language}
Presburger's language is a language with variables running over  $\Zgordon$,
and symbols '$+$', '$\le$', '$0$', '$1$', and for each $d=2,3,4,\dots$, a symbol '$\equiv_d$' to denote the binary 
relation $x\equiv y \pmod{d}$, together with all the symbols for quantifiers,
logical operations and parentheses, as above.
Note the absence of the symbol for multiplication. 

Since multiplication is not allowed, definable functions have to be linear 
combinations of piecewise-linear and periodic functions (where the period is a 
vector in $\Zgordon^n$, and $n$ is the number of variables).

\subsection{The language of Denef-Pas}
The language of Denef-Pas is designed for valued fields. It is a \emph{three-sorted language}, meaning that
it has three sorts of variables. Variables of the first sort run over the valued field, variables of the second 
sort run over the value group (for simplicity, we assume that the value group is $\Zgordon$), and variables of the
third sort run over the residue field. 

The symbols for this language consist of the symbols of the language of rings for the residue field sort,
Presburger's language for the $\Zgordon$-sort, and the language of rings for the valued field sort, together
with two additional symbols:
$\ord(x)$ to denote a function from the valued field sort to the $\Zgordon$-sort,  and $\ac(x)$ to denote a function from
the valued field sort to the residue field sort.
These functions are called the {\bf valuation map}, and the 
{\bf angular component map}, respectively.  
We also need to add the symbol {`}$\infty$' to the value sort, to denote the valuation of {`}$0$'
(so that {`}$\ord(0)=\infty$' has the {`}true' value in every structure). 

A valued field $K$ \emph{together with the choice of the uniformizer of the valuation on $K$} is a 
structure for Denef-Pas language. 
In order to match the formulas in Pas's language with their interpretations in 
its structure $K$, we need to give a meaning to the symbols {`}$\ord$' and 
{`}$\ac$' in the language. 

The function $\ord(x)$ stands for the valuation of $x$. 
 In order to provide the interpretation for the symbol {`}$\ac(x)$', 
we have to   
fix a uniformizing parameter $\varpi$. 
The valuation on $K$ is normalized 
so that
$\valgordon(\varpi)=1$.
If $x\in\rigordon_{K}^{\ast}$ is a unit, there is a natural definition
of $\ac(x)$ -- it is the reduction of $x$ modulo the ideal $(\varpi)$.
Define, for $x\neq 0$ in $K$, $\ac(x)=\ac(\varpi^{-\ord(x)}x)$,
 and $\ac(0)=|0|=0$.

For convenience, a symbol for every rational number is added to the valued field sort, so that we could have formulas with coefficients in $\Qgordon$.

Sometimes, when the category of fields under consideration is restricted to 
all fields containing a fixed ground field $k$, one can add a symbol for each 
element of $k((t))$ to the valued field sort. This enlarges the class of definable sets. 
In order to make distinctions between various settings, we will explicitly 
talk  of ``formulas with coefficients in $k((t))$ (or in $k[[t]]$)'' in such 
cases. Note that in any case, for an arbitrary  field $K$ containing $k$, coefficients from $K$, or $K((t))$, are not allowed (otherwise this would have been meaningless -- we want to 
consider the sets of points satisfying a given formula for the varying 
fields $K$). Given a local field $K$ containing $k$ with a uniformizer 
$\varpi$, one can make a map from $k((t))$ to $K$ where $t\mapsto\varpi$ 
(this will be discussed in detail in Section~\ref{sec:back}). In this sense, 
$t$ plays the role of the uniformizer of the valuation, to some extent.

We will talk in detail about interpreting formulas in different structures
in Section~\ref{sec:back}.

\subsection{Definable subassignments}
Here we introduce the terminology that conveniently puts 
the set of points defined by an 
interpretation  of a logical formula over a given field on the same 
footing with, say, the set of points of an affine variety. To do that, we use the 
language of functors.

We fix a ground field $k$ of characteristic $0$. 
For most applications, one can think that $k=\Qgordon$. 
Denote by $\text{Field}_k$ the category of fields containing $k$.
Any variety $X$ over $k$ defines a functor -- its functor of points -- 
from $\text{Field}_k$ to 
the category of sets, by sending every field $K$ containing $k$ to $X(K)$.
This functor will be denoted by $h_X$.

\begin{definition} We will denote by  $h[m,n,r]$ 
(or $h_{\A^m_{k((t))}\times \A^n_k\times \Zgordon^r}$)
\footnote{
Even though the objects whose volumes we would like to compute correspond to
subsets of affine spaces over the \emph{valued field}, it is
very useful to
have a formalism that allows us to deal with valued-field, residue-field, and 
integer-valued variables at the same time.
 One of the advantages of doing that
is being able to look at definable families with integer-valued or 
residue-field valued parameters. This is the reason that this functor 
plays a fundamental role.}  
the functor from the category 
$\text{Field}_k$ to the category of sets defined by 
$$
 h_{\A^m_{k((t))}\times \A^n_k\times \Zgordon^r}(K)=K((t))^m\times K^n\times \Zgordon^r.
$$ 
\end{definition}

For example, $h[1,0,0]$ is the functor of points of $\A^1_{k((t))}$, and 
$h[0,0,0]$ is a functor that assigns to each field $K$  
a one-point set. We will usually write $h_{\spec k}$ for $h[0,0,0]$. 

\begin{definition} 
Let $F:{\mathcal C}\to \underline{\text{Sets}}$ be a functor from a 
category ${\mathcal C}$ to the category of sets.
A {\bf subassignment} $h$ of $F$ is a collection of 
subsets $h(C)\subset F(C)$, one for each 
object $C$ of ${\mathcal C}$. 
\end{definition}
Note that a subassignment does not have to be a subfunctor 
(that is, we are making no requirement that a morphism 
between two objects $C_1$ and $C_2$ in ${\mathcal C}$ has to correspond to 
a map between the corresponding sets $h(C_1)$ and $h(C_2)$). 

The subassignments will replace formulas in the same way that 
functors can replace varieties. When we talk about formulas, we will mean 
logical formulas built using {the} Denef-Pas language (so in particular, {we use} the  
language of rings for the residue field, and Presburger language for $\Zgordon$). 

\begin{definition}
A subassignment $h$ 
of $h[m,n,r]$ is called {\bf definable} if there exists a formula $\phi$ 
in the language of Denef-Pas with coefficients in $k((t))$, with $m$ free
variables of the valued field sort, with coefficients in $k$ and $n$ free variables of the residue field sort, and $r$ free variables of the value sort, such 
that for every $K$ in $\text{Field}_k$, $h(K)$ is the set of all points in 
$K((t))^m\times K^n\times \Zgordon^r$ satisfying $\phi$. 
\end{definition}

\begin{definition} A {\bf morphism of definable subassignments} $h_1$ and $h_2$
is { a definable subassignment $F$ such that $F(C)$ is the graph of a
function from $h_1(C)$ to $h_2(C)$ for each 
object $C$.} 
The {\bf category of definable subassignments} of some $h[m,n,r]$ 
will be denoted {\bf $\de_k$}. 
\end{definition}

\subsubsection{Relative situation}
If $S$ is 
an object in $\de_k$, one can consider the category of definable 
subassignments equipped with a morphism to $S$, denoted by $\de_S$ 
(the morphisms being the morphisms over $S$). 
More precisely, we could say that the objects 
are morphisms $[Y\to S]$ with $Y\in\de_k$, and morphisms are 
commutative triangles
\begin{equation*}
\xymatrix{
&W \ar[d] \ar[r] & Y \ar[ld] \\
 & S& \quad.
}
\end{equation*}

We denote by $S[m,n,r]$ the subassignment
$$
S[m,n,r]:=S\times h_{\A^m_{k((t))}\times \A^n_k\times \Zgordon^r};$$
This is an object of $\de_S$, the morphism to $S$ being the projection onto the 
first factor.

Finally, for $S$ an object in $\de_k$, there is the category 
of {\bf R-definable subassignments over $S$}, denoted by 
{\bf $\rde_S$} ($R$ stands for ``residue''). 
The objects of $\rde_S$ are definable subassignments
of $S[0,n,0]$ for some integer $n\ge 0$ (with a morphism to $S$ coming from the projection onto the first 
factor), and morphisms are morphisms over $S$.
Note that this abbreviation says that the objects in $\rde_S$ can have extra variables of the residue field sort, but \emph{no extra variables of the valued field sort nor the value group sort}, compared to $S$ itself.

\begin{example}
{\bf The category $\rde_{\spec k}$.}
By definition, the category $\rde_{h_{\spec k}}$ consists of definable subassignments with variables ranging only over the residue field (and therefore definable in the language of rings).
Note that if the formulas defining the subassignments in $\rde_{\spec k}$
had been 
quantifier-free, then they would essentially define 
{constructible sets} 
over $k$. 
Depending on $k$, since quantifiers are allowed, this category may be richer, 
but in many cases 
there is a map from it to a category of geometric objects over the residue 
field, as discussed in Section~\ref{sec:back}. 
\end{example}

The category $\rde_{\spec k}$ (and more generally, $\rde_S$ where $S$ 
is a definable subassignment) is going to play a very important role 
in the theory.
In the next section, we will associate with each definable 
subassignment its motivic volume that will be, essentially, an element of 
the {Grothendieck} ring (defined in the next section) of the category 
$\rde_{\spec k}$.

\subsubsection{Points on subassignments, and functions}
By definition, a point on a definable subassignment $Y\in \de_k$ is a 
pair $(y_0, K)$ where $K\in \text{Field}_k$, and $y_0\in Y(K)$.

Given any definable morphism
$\alpha: S\to Z$, where both $S$ and $Z$ are definable {subassignments}, 
there is a corresponding function from the set of 
points of $S$ to the set of points of $Z$. The function and the morphism define each other uniquely, so we can identify them.
In the special case $Z=h[0,0,1]$, the resulting function is integer-valued, so we will say that such a morphism  
is an integer-valued function on the subassignment $S$.

\subsection{Grothendieck rings}\label{sub:gr.rings}
There are several Grothendieck rings used in various {constructions} of motivic measure. The first one is the Grothendieck ring of the category of varieties 
over $k$, $K_0(\var_k)$. Its elements are formal linear combinations
with coefficients in $\Zgordon$ of isomorphism classes of varieties (with formal 
addition) modulo the 
natural relation $[X\setminus Y]+[Y]=[X]$;
the product operation comes from the product in the 
category $\var_k$.

Another Grothendieck ring that is sometimes used is $K_0(\mot_k)$ -- the Grothendieck ring of the category of Chow motives over $k$. 
(We will not talk about Chow motives here, see \cite{Scholl} for an 
introduction).
This is the ring constructed in the same way, but from the category of Chow motives rather than varieties over $k$.

These rings have an element (corresponding to the class of the affine line)
that plays a special role in the theory of motivic integration. It is always 
denoted by $\lef$. The notation comes from Chow motives, where $\lef$ stands 
for the so-called Lefschetz motive $\lef=[{\mathbb P}^1]-[pt]$ (see 
\cite{Scholl}).
In $K_0(\var_k)$, $\lef$ stands for $[\A^1]$. 
It is a difficult theorem (Gillet and Soul\'e, \cite{gillet-soule}, 
and Guill\'en and Navarro Aznar) 
that there exists a natural map from $K_0(\var_k)$ to $K_0(\mot_k)$.
\footnote{The meaning of ``natural'' here is the following. Chow motives are,
formally, equivalence classes of triples $(X, p, n)$, where $X$ is a variety, 
$p$ is an idempotent correspondence on $X$ (one can think of it as a projector from $X$ to itself), and $n$ is an integer. Every \emph{smooth projective} variety $X$
naturally corresponds to the Chow motive $(X, \text{id}, 0)$. The content of the theorem is to extend this map to the elements of $K_0(\var_k)$ that are not necessarily linear combinations of isomorphism classes of smooth projective varieties.}  
Under this map, the class of the affine line corresponds to $\lef$ 
(see \cite{Scholl}), thus justifying the notation.  
The image of this map will be denoted by $K_0^{mot}(\var_k)$, and it will play an important role in Section~\ref{sec:back}.

One can also make Grothendieck rings of the categories of subassignments that we have considered above. Note that one can define set-theoretic operations on subassignments in a natural way, \eg, 
$(h_1\cup h_2)(K):=h_1(K)\cup h_2(K)$, {\it etc}.
Let $S$ be a definable subassignment. 
One can make the ring $K_0(\rde_S)$: its elements are formal linear 
combinations of isomorphism classes of objects of $\rde_S$, modulo the relations $[(Y\cap X)\to S]+[(Y\cup X)\to S]=[Y\to S]+[X\to S]$, and 
$[\emptyset\to S]=0$. With the natural operation of addition, $K_0(\rde_S)$ is an abelian group; cartesian product gives it a structure of a ring. 

\begin{remark}
Note that when making a Grothendieck ring, we first replace the objects
of a category by equivalence classes of objects. By changing the notion of equivalence (for example, making it more crude), one can define the rings where various important invariants take values. We shall see in Section~\ref{sub:pseudo} that in order to get a version of motivic integration that specializes to $p$-adic 
integration, we need to replace equivalence by \emph{equivalence on pseudofinite fields}.
\end{remark}

\subsubsection{Dimension} 
Before we can talk about measure theory for objects of $\de_k$, we need a dimension theory.  
Recall that each subassignment has valued-field, residue-field, and value-group variables. The notion of dimension takes into account only the valued-field 
variables (this {fits well} with the measure theory we are about to describe 
 since the measure on $K^n\times \Zgordon^r$ is 
going to be {essentially} the counting measure, as {we will} see below).

First, note that each subvariety $Z$ of $\A^m_{k((t))}$ naturally gives a 
subassignment $h_Z$ of $h[m,0,0]$ by $h_Z(K):=Z(K((t)))$.  
For $S$ a subassignment of $h[m,0,0]$, we define {the} {\bf Zariski closure of $S$} to be 
the intersection $W$ of all subvarieties $Z$ of $\A^m_{k((t))}$ such that $h_Z$ contains $S$. Then {the} dimension of $S$ is defined to be the dimension of $W$.

In general, if $S$ is a subassignment of $h[m,n,r]$, the dimension of $S$ is defined to be the dimension of the projection of $S$ onto the first component
$h[m,0,0]$.

\begin{proposition}\cite[Prop.~3.4]{CL.expo}
Two isomorphic objects of $\de_k$ have the same dimension.
\end{proposition} 

Note that {definable} subassignments are closely related to analytic manifolds.
See \cite[\S~3.2]{CL} for a detailed discussion.

\subsection{Constructible motivic Functions}\label{sub:constr.f}
\subsubsection{The ring of values}\label{subsub:ring_of_values}
Let $\lef$ be a formal symbol (later it will be associated with the class of the affine line in an appropriate Grothendieck ring). In Section~\ref{sec:back}, 
it will be matched with $q$ -- the cardinality of the residue 
field -- but for now it is just a formal symbol.
 
{Consider} the ring
$$
A=\Zgordon\left[\lef, \lef^{-1}, \left(\frac{1}{1-\lef^{-i}}\right)_{i>0}\right].
$$ 
For all real $q>1$, there is a homomorphism of rings
$\nu_q:A\to \R$ defined by $\lef\mapsto q$.
Note that if $q$ is transcendental, then $\nu_q$ is injective.

This family of homomorphisms gives a partial ordering on $A$: for $a, b\in A$,
set $a\ge b$ if for every real $q\ge 1$ we have $\nu_q(a)\ge \nu_q(b)$.
Note that with this ordering, $\lef^i$, $\lef^i-\lef^j$ with $i>j$, and 
$\frac{1}{1-\lef^{-i}}$ with $i>0$ are all positive, but for example, $\lef-2$ is not positive.

\subsubsection{Constructible motivic functions}
In the $p$-adic setting, the smallest class of functions that one would definitely like to be able to integrate is built from two kinds of functions: characteristic functions of measurable sets, and functions 
of the form $q^{\alpha}$, where $q$ is the cardinality of the residue field, 
and $\alpha$ is a characteristic function of a measurable set (these appear as absolute values of the functions of the first kind).
Keeping this in mind, let us define constructible motivic functions.

Let $S\in \de_k$ be a definable subassignment.
The ring of constructible motivic functions {on $S$} is built from 
two basic kinds of functions.

The first kind are definable functions with values in $\Zgordon$, and functions 
of the form $\lef^{\alpha}$, where $\alpha$ is a definable function on $S$ with values in $\Zgordon$
(these functions can be thought of as functions with values in $A$). In particular, 
this collection of functions includes characteristic functions of
 definable subsets of $S$. Let us denote the ring of $A$-valued 
functions on $S$ 
generated by functions of these two kinds, by ${\mathcal P}(S)$. 

The second kind 
of definable functions on $S$ do not look like functions at all, 
at the first glance.
Formally, they are the elements of the Grothendieck ring $K_0(\rde_S)$, as defined in Section~\ref{sub:gr.rings}.
However, if we think {of specialization} to $p$-adic integration, we see that
once we have fixed a local field $K$ with a (finite) residue field $\F_q$,  
an element of $[Y\to S]\in \rde_S$ gives an integer-valued function on $S$ by 
assigning to each point
on $x\in S(K)$ the cardinality of the {fibre} of $Y$ over $x$. 
Note that the {fibre} of $Y$ over $x$ is a subset of $\F_q^n$ for some $n$;  
in particular, it is finite.

The reason these functions need to be included from the very beginning is that the 
motivic integral will take values in a ring containing $K_0(\rde_S)$, and we need to be able to integrate a  function of two variables with respect to one of the variables, and get a function of the 
remaining variable that is again integrable.

To put together the two kinds of functions described above, 
note that characteristic functions of definable subsets of $S$ naturally 
correspond to elements of $\rde_S$: ${\bf 1}_Y$ corresponds to 
$[Y\to S]\in \rde_S$.
Let ${\mathcal P}^0(S)$ be the subring of ${\mathcal P}(S)$ generated
by the constant 
function $\lef_S-1_S$ (where $\lef_S=[S\times \A_{k((t))}^1\to S]$, and 
$1_S=[S\times h_{\spec k}\to S]$), and the functions of the 
form 
${\bf 1}_Y$, where $Y$ is a definable subassignment of $S$.
We can form the tensor product of the ring ${\mathcal P}(S)$ and the ring $K_0(\rde_S)$:
$$
{\mathcal C}(S):={\mathcal P}(S)\otimes_{{\mathcal P}^0(S)}K_0(\rde_S).
$$
This is the ring of {\bf constructible motivic functions on $S$}.
We refer to 
\cite[\S~3.2]{CL.expo} for details.

Finally, one defines {\bf constructible motivic Functions on $S$} as equivalence classes of
elements of ${\mathcal C}(S)$ ``modulo support of smaller dimension''.
See \cite[\S~3.3]{CL.expo} for a precise definition and discussion why this needs to be done. We will think of constructible motivic Functions as
functions defined almost everywhere (which is quite reasonable in the context of any integration theory).

\subsection{Summary}
Let $k$ be the base field, {\eg}, $k=\Qgordon$.
To summarize, 
instead of measurable sets we have definable {subassignments}; instead of functions -- constructible motivic Functions; and instead of numbers as values of the measure -- elements of a suitable Grothendieck ring (either of varieties, or of 
Chow motives, or of $\rde_k$, depending on the context).

The measure theory and its relation to $p$-adic measure 
is summarized by the diagram.
 
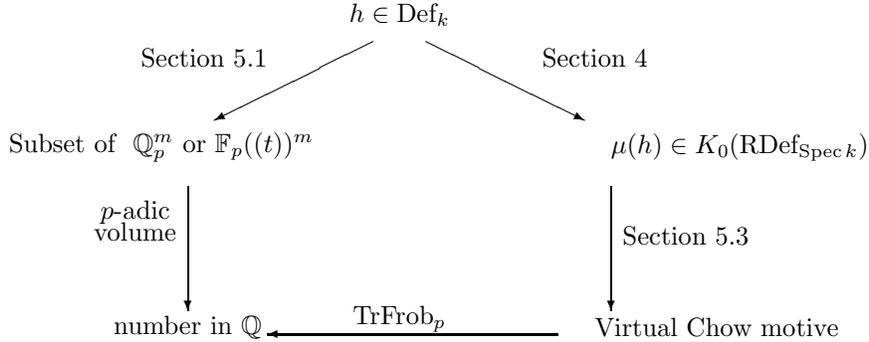
\begin{picture}(400,160)(50,30) 
\put(110,210){\makebox(0,0)} 
\put(190,150){\vector(-2,-1){60}} 
\put(150,144){\makebox(0,0)[r]{Section~\ref{sub:interpr}}} 
\put(200,160){\makebox(0,0) 
{$h\in \de_k$}} 
\put(210,150){\vector(2,-1){60}} 
\put(254,144)
{\makebox(0,0)[l]{Section~\ref{sec:CL}}} 
\put(280,110) 
{\makebox(0,0)[l]{$\mu(h)\in K_0(\rde_{\spec k})$}} 
\put(120,95){\vector(0,-1){47}} 
\put(110,110){\makebox(0,0){\text{Subset of } $\Qgordon_p^m$ or $\F_p((t))^m$}} 
\put(280,95){\vector(0,-1){47}} 
\put(120,42){\makebox(0,0){number in $\Qgordon$}} 
\put(100,90){\makebox(0,0)[t] {$p$-adic}}\put(100,78){\makebox(0,0){volume}} 
\put(320,42){\makebox(0,0){Virtual Chow motive}} 
\put(260,39){\vector(-1,0){110}} 
\put(200,50){\makebox(0,0)[t]{$\text{TrFrob}_p$}} 
\put(305,80){\makebox(0,0)[t] {\ \ Section~\ref{sub:comp}}} 
\end{picture}

We describe the arrow from subassignments to elements of $K_0(\rde_k)$ 
in the next section (this is what motivic integration developed in \cite{CL} essentially amounts to). We explain the relationship with $p$-adic integration in 
Section~\ref{sec:back}, as indicated on the diagram.

\begin{remark}
As we will see, for the sets that come from definable subassignments,
the value of the $p$-adic measure, that is claimed to be 
in $\Qgordon$ (in the bottom left corner of this diagram) in fact lies in
$\Zgordon\left[\frac1p, \left(\frac{1}{1-p^{-i}}\right)_{i>0}\right]$.
\end{remark}

In this diagram, one can make a choice for the collection of fields that 
appears in the upper left-hand corner.  One natural collection of local fields
 would be the collection ${\mathcal A}_F$ of all finite degree field extensions of all  non-archimedean 
completions of a given global field $F$ 
(in that case, one adds to Denef-Pas
language constant symbols for all elements of $F$). Another natural collection
is the collection of all function fields $\F_q((t))$. One of the most
spectacular applications of motivic integration is the 
\emph{Transfer Principle}
that allows to transfer identities between these two collections of fields.
We talk more about this in Section~\ref{sec:applic}.

\section{Motivic integration as pushforward}\label{sec:CL}
The main difference between motivic integration developed by
Cluckers and Loeser \cite{CL} and the older theories is that in \cite{CL} integration, by definition, is \emph{pushforward} 
of morphisms, in agreement with Grothendieck's philosophy.

Let $f:S\to W$ be a morphism of definable subassignments.
We have described the rings of constructible motivic functions 
${\mathcal C}(S)$ and ${\mathcal C}(W)$ on $S$ and $W$, respectively. The goal
is to define a morphism of rings $f_{!}:\calf(S)\to \calf(W)$ that 
corresponds
to integration along the fibres of $f$.\footnote{In reality, the situation is more complicated because, naturally, not all constructible functions are integrable. Accordingly, one needs to define a class of integrable functions. We say a few words about it in Section~\ref{sub:carpet}, but for now we will ignore this issue 
for simplicity of exposition.} 

To make the situation more manageable, the operation of pushforward is defined for various types of projections and injections, keeping in mind that  a 
general morphism can be represented as the composition of a projection and an
 injection by considering its graph.

Naturally, pushforward for injections is extension by zero, and the 
interesting 
part is the projections.
There are three kinds of projections: forgetting the valued-field, residue-field, or $\Zgordon$-valued variables.
It is a nontrivial proposition that the three kinds of variables are independent, in the sense that you can pushforward along these projections in any order. 

In order to understand the theory completely, one needs to read \cite{CL}. 
Here we only aspire to sketch integration with respect to one valued-field variable. The idea is to break up the domain of integration into simpler sets (the cells), and define the integral on each cell. Then one can repeat this procedure inductively to integrate along all the variables and get the volume. The hardest part of the theory is a collection of the statements of Fubini type that allow to permute the order of integrals with respect to the valued-field variables.

Throughout this section, we fix the ground field $k$ and let $S\in\de_k$ be a 
definable subassignment (of some $h[m,n,r]$). 
We start with the exposition of cell decomposition theorem, which is the main 
tool of the construction.

\subsection{The Cell Decomposition Theorem}
Cell decomposition theorem is a very powerful theorem with many striking 
applications. 
The  article \cite{Denef} gives a beautiful exposition of 
$p$-adic cell decomposition (with a slightly 
more restrictive definition of cells) and its 
applications to questions about rationality of Poincar\'e series.
Here we will focus, instead, on examples illustrating 
the technical side of the cell decomposition 
used in the construction of the motivic measure.
 
Before we state the theorem, let us consider a simple example 
of a $p$-adic integral. 

\subsubsection{A motivating example}\label{subsub:example}
Consider the integral, depending on a parameter $x\in \Zgordon_p$: 
$$\int_{\Zgordon_p}|t^3-x||\dt|.$$
Let us calculate this integral by brute force, as a computer could have 
done it. We assume that $p>3$. 

First, consider the easiest case $x=0$.
Then the domain breaks up into infinitely many ``annuli'' $A_i$ 
on which the function 
$|t^3|$ is constant. (Even though each one of the sets $A_i$ lives on the 
line, we call it an annulus because it is a difference of two $1$-dimensional balls of radii $p^{-i}$ and $p^{-(i+1)}$ respectively). 

 The volume of each annulus is:
\begin{equation*}
\begin{aligned}\mu(A_i)&=\mu(\{t\mid |t|=p^{-i}\})=\mu(\{t\mid |t|\le p^{-i}\})-
\mu(\{t\mid |t|\le p^{-(i+1)}\})\\
&=p^{-i}-p^{-(i+1)}=p^{-(i+1)}(p-1).
\end{aligned}
\end{equation*}
Then the value of the integral for $x=0$ is the sum of the geometric series:
$$
\int_{\Zgordon_p}|t^3||\dt|=\sum_{i=0}^{\infty}p^{-3i}p^{-(i+1)}(p-1)=
\frac{p-1}p\frac{1}{1-p^{-4}}.
$$

Now let us turn to the case of general $x$. 
If $\ord(x)$ is not divisible by $3$, then  for any $t$, we have 
$|t^3-x|=\max(|t^3|, |x|)$, and so 
the domain of integration breaks up into two parts: the part where $|t^3|$ 
dominates, and the part where $|x|$ dominates. 
The integral over each part is easily 
reduced to the sum of a geometric series,  and 
we omit the details.

The most interesting case is the case where $\ord(x)=3k$ for some integer $k$: 
in this case, along with the two ``easy'' integrals similar to the previous case (which we omit) there is also the integral  over the set 
$B=\{t\mid |t^3|=|x|\}$.   
This case breaks up further into three subcases:
\begin{enumerate}
\item $x$ is not a cube;
\item $x$ is a cube, and there is one cube root of $x$ in $\Zgordon_p$;
\item $x$ is a cube, and there are three cube roots. 
\end{enumerate}
 
Case (1) is also easy to finish, because in this case the formula 
$|t^3-x|=\max(|t^3|, |x|)$ still holds for all $t$.
We will focus on the cases (2) and (3), which are the most interesting.
If '$\exists y \mid x=y^3$' holds, then the number of solutions to this 
equation in $\Zgordon_p$ depends on $p$: for example, 
there is only one root in $\Zgordon_5$, and three roots in $\Zgordon_7$. 
Let us consider the case with $3$ roots first.

We can write $t^3-x=(t-y_1)(t-y_2)(t-y_3)$.
Suppose $t\in B$. 
First, consider the subset $B_0$ of $B$ that consists of the points $t$ 
such that $\ac(t)\neq \ac(y_i)$, $i=1,2,3$. On this set, $|t-y_i|=p^{-k}$, 
and 
\begin{equation*}
%\label{eq:B0}
\int_{B_0}|t^3-x|=p^{-3k}\mu(B_0)=p^{-3k}((p^{-3k}-p^{-3(k+1)})-3p^{-3(k+1)})=
p^{-6k}-4p^{-(6k+1)}.
\end{equation*}
Finally, consider the sets $B_i=\{\ac(t)=\ac(y_i)\}$, $i=1,2,3$. 
It is enough to understand the integral over one of them, say, $B_1$.

The set $B_1$ is defined by 
$$
B_1=\{t\mid \ord(t)=k=\ord(b_1) \wedge \ac(t)=\ac(y_1)\}=
\{t\mid \ord(t-y_1)\ge (k+1)\}.$$
The integral over $B_1$ becomes an infinite sum (indexed by the degree of congruence between $t$ and $y_1$, which we denote by $m=\ord(t-y_1)$):
\begin{align}\label{eq:B_1}
\int_{B_1}|t^3-x| &=\sum_{m=k+1}^{\infty}p^{-m}p^{-2k}(p^{-m}-p^{-(m+1)})\\
&=(1-p^{-1})p^{-2k}p^{-2(k+1)}(1-p^{-2})^{-1}.
\end{align}
 
From here it is easy to get the final answer, and easy to do the case of one cube root of $1$ in the field. The main point here is that in each case, 
the integral boils 
down to a few geometric series with a power of $p$ as the ratio, and a few finite sums. As we will see, this is a very general pattern.

The interesting part of the final answer for the case when the 
parameter $x$ is a cube, $\ord(x)=3k$: 
$$
\int\limits_{\{t \,\mid\, 3\ord(t)=\ord(x)\}}\hspace*{-2.5em}|t^3-x||\dt|=
\begin{cases}
3(1-p^{-1})\displaystyle\frac{p^{-4k-2}}{1-p^{-2}}+p^{-6k}-4p^{-(6k+1)}, &
p\equiv 1 \pmod{3};\\
(1-p^{-1})\displaystyle\frac{p^{-4k-2}}{1-p^{-2}}+p^{-6k}-2p^{-(6k+1)}, &
p\equiv 2 \pmod{3}.
\end{cases}
$$

\subsubsection{The definition of cells}
The general idea behind cell decomposition is to present every definable 
set as a fibration over some definable set of dimension one less
(called the basis) with fibres that are $1$-dimensional $p$-adic balls.

\begin{definition} 
Let $S$ be a definable subassignment.
Let $C\subset S$ be a definable subassignment of $S$, and 
let $c:C\to h[1,0,0]$, $\alpha:C\to \Zgordon$, $\xi:C\to h_{\Ggordon_m, k}$ 
be definable morphisms.  
Denote by  {\bf $Z_{C, \alpha,\xi, c}$} a subassignment of $S[1,0,0]$ 
defined by $y\in C$, $\ord(z-c(y))=\alpha(y)$, $\ac(z-c(y))=\xi(y)$. 
The subassignment  $Z_{C,\alpha, \xi, c}$ is a basic $1$-cell.
We will refer to the subassignment $C$ as its {\bf basis} and to the function 
$c$ as its {\bf centre}.  
\end{definition}

When doing cell decomposition, we will also need to be able to have some 
{pieces of smaller dimension}. This is the {idea} behind the next 
definition. 

\begin{definition}
In the context of the previous definition,
denote by  
{\bf $Z_{C,c}$}  the subassignment of $S[1,0,0]$ defined by the 
formula $y\in C, z=c(y)$. This is a {basic $0$-cell}.
This is a subassignment of the same dimension as $C$; essentially, 
it is a copy of $C$ that sits in a space of dimension {one greater}. 
\end{definition}

These basic cells are simple enough to work with, but not yet versatile enough
for cell decomposition to work.  
We need to modify the definition of cells by allowing extra 
residue field and integer-valued variables, and letting the points of the 
cell live on different ``levels'' according to the values of these variables. 

\begin{definition}\label{def:cells}
Let $S$ be a definable subassignment, let $s,r$ be some non-negative integers,
and let $\pi$ be the projection 
$\pi:S[1,s,r]\to S[1,0,0]$ onto the first factor.
A definable subassignment $Z\subset S[1,0,0]$ is called a {\bf $1$-cell} if
there exists an isomorphism of definable subassignments 
(called a {\bf presentation})
$$\lambda:Z\to Z_{C,\alpha, \xi, c}\subset S[1,s,r]$$ 
for some $s, r \ge 0$, some basis $C\subset S[0,s,r]$, 
such that $\pi\circ \lambda$ is the identity on $Z$.
\[
\xymatrix{
Z_{C,\alpha,\xi,c} \ar[rr]^{\hookrightarrow}&& \ar[d]^{\pi} S[1,s,r]\\
& \ar[ul]^{\lambda} Z \ar[r]^{\hookrightarrow} & \ar[d] S[1,0,0]\\
C \ar[r]^{\hookrightarrow} & S[0,s,r] \ar[r] & S
}
\]

\end{definition}

A similar definition applies to $0$-cells, with the only change that 
the isomorphism $\lambda$ is between $Z\subset S[1,0,0]$ and a $0$-cell
$Z_{C,c}\subset S[1,s,0]$ with some basis $C\subset S[0,s,0]$ (in particular, no extra $\Zgordon$-valued variables allowed).

%The only problem with those diagrams is that the way the definitions are 
%set up, the $Z's$ are really subsets, not just abstract objects that  have an embedding into, the corresponding sets $S[...]$. Because of this, the prettier diagrams that you made are a little misleading... I agree though that my diagrams don't look good. But I left them in for now.

\begin{example} Take $S=\spec k$.
We can write the line $h[1,0,0]$ as the union of a 
$0$-cell $h_{\spec k}$ and a 
$1$-cell $Z=\A^1_{k((t))}\setminus \{0\}$ (this is not a very precise 
notation for a subassignment but this makes the meaning more clear).   
Let us see precisely why $Z$ is indeed a $1$-cell. Let us define the subassignment $Z_{C,\alpha,\xi, c}$ and the presentation $\lambda$ 
required by the definition. We have the freedom to choose the number of extra 
residue field 
and $\Zgordon$-valued variables to introduce. Let us make 
$Z_{C,\alpha,\xi, c}$ a subassignment of $h[1,1,1]$.
As the basis, we take the subassignment  $C$ of $h[0,1,1]$ defined by 
$\eta\neq 0$ (recall that $h[0,1,1]$ stands for $\A^1_k\times \Zgordon$). 
We call the residue field variable $\eta$, and the $\Zgordon$-variable $r$. Let
$c(\eta,r)=0$ be the constant zero function from  
$h[0,1,1]$ to $h[1,0,0]$ (\ie, to $\A_{k((t))}^1$), 
and let $\xi(\eta, r)=\eta$, $\alpha(\eta, r)=r$.
Now let $Z_{C,\alpha,\xi, c}$ be the subassignment of
$h[1,1,1]$ (denote the variables by $(z, \eta, r)$) 
defined by $\ord(z)=r$, $\ac(z)=\eta$.

The presentation $\lambda:Z\to Z_{C,\alpha,\xi, c}$ is given by
$\lambda(z)=(z,\ac(z),\ord(z))$.
The projection $\pi$ is the projection onto the first factor from 
$h[1,1,1]$ to $h[1,0,0]$ (that is, we forget the extra residue field and $\Zgordon$-variables). Clearly, $\pi\circ\lambda$ is the identity on $Z$. 

One way to think about it is to imagine that we have placed different 
points on the affine line (without $0$) over  the valued field on 
different ``layers'' 
indexed by their valuations and angular components. 
\end{example}

\subsubsection{Cell Decomposition Theorem}

\begin{theorem}[\cite{CL}, Th.~7.2.1]
Let $X$ be a definable subassignment of $S[1,0,0]$ with $S$ in $\de_k$.
%Then
\begin{enumerate}
\item The subassignment $X$ is a finite disjoint union of cells.
\item For every constructible function $\varphi$ on $X$ there exists 
a finite partition of $X$ into cells $Z_i$ with presentations 
$(\lambda_i, Z_{C_i})$ such that $\varphi\mid_{Z_i}$ is the pullback by 
$p_i\circ \lambda_i$ of a constructible function $\psi_i$ on $C_i$, where
$p_i$ is the projection $p_i:Z_{C_i}\to C_i$.
This is called the {\bf cell decomposition adapted to $\varphi$}. 
\end{enumerate}
\begin{center}
\[
\xymatrix{
\ar[dd]^{p_i} Z_{C_i,\alpha_i,\xi_i,c_i} \ar[rr]^{\hookrightarrow}&& \ar[d] S[1,s_i,r_i]\\
& \ar[ul]^{\lambda_i} Z_i \ar[r]^{\hookrightarrow} & \ar[d] S[1,0,0]\\
C_i \ar[r]^{\hookrightarrow} & S[0,s_i,r_i] \ar[r] & S
}
\]
\end{center}
\end{theorem}

\begin{example}\label{ex:cells}
Let us consider the cell decomposition adapted to the function 
$\varphi(x,t)=|t^3-x|$ with respect to the $t$-variable
(see our ``motivating example'',  
Section~\ref{subsub:example}).
Note that $|t^3-x|_p=p^{-\ord(t^3-x)}$, so it is natural to define the
corresponding
$A$-valued function (which we will also denote by $\varphi(x,t)$ in 
this example) by  $\varphi(x,t)=\lef^{-\ord(t^3-x)}$. (The details about the interpretation of constructible motivic functions  will appear below in 
Section \ref{sec:back}.)

As in Section~\ref{subsub:example}, it is convenient to consider the case
$x=0$ separately.
In our language, $\phi(x,t)$ is a function on $h[2,0,0]$.
We split the domain into the two subassignments defined by $x\neq 0$ and $x=0$.
We only deal with  the part $x\neq 0$ as it is more interesting.
 
First, consider the subassignments $h_1$ and $h_2$  defined by 
$3\ord(t)<\ord(x)$ and $3\ord(t)>\ord(x)$, respectively.
  
On $h_2$ we have  $f(x,t)=|x|$. Since $f(x,t)$ is independent of $t$, this is the easiest part:
$h_2$ is a single cell and the function $\psi$ is 
$\lef^{-\ord(x)}$. The details of the presentation are 
left to the reader.

The subassignment $h_1$ is a single cell as well. Indeed, on $h_1$, we have 
$f(x,t)=|t^3|$.  To define the basis  $C$, we add extra value sort variables 
for $\ord(x)$ and $\ord(t)$, and an extra residue field variable for $\ac(t)$: 
formally, let $C$ be the 
subassignment of $h[1,1,2]$ defined by the formula 
$$\phi(x, \eta, \gamma_1,\gamma_2)=\text{`} (x\neq 0) 
\wedge (\gamma_1=\ord(x))\wedge (3\gamma_2<\gamma_1)
\text{'}.
$$

Let the centre $c: C\to h[1,0,0]$ be the zero function, let 
$\alpha:C\to \Zgordon$ be the function 
$(x,\eta, \gamma_1,\gamma_2)\mapsto \gamma_2$,
and let $\xi(x,\eta,\gamma_1,\gamma_2)=\eta$ (so that $\xi$ is a function from 
$C$ to $\Ggordon_m$).
Let $Z_{C,\alpha,\xi, c}$ be the subassignment of 
$h[2,1,2]$ defined by 
$$\phi_1(x,t,\eta,\gamma_1, \gamma_2)=\text{`}
(\ord(t)=\gamma_2) \wedge (\ac(t)=\eta)
\text{'}.
$$
The presentation $\lambda: h_1\to Z_{C,\alpha,\xi, c}$ is given by
$$\lambda(x,t)=(x,t,\ac(t),\ord(x),\ord(t)).
$$ 

Finally, let $\psi$ be the function on $C$ (with values in the ring $A$ 
of Section~\ref{subsub:ring_of_values}) defined by 
$\psi(x,\eta,\gamma_1,\gamma_2)=\lef^{-\gamma_1}$.

Then, clearly, on $h_1$ 
our function $\lef^{-\ord(t^3-x)}$ is the pullback of $\psi$ by
$p\circ\lambda$.

Now let us consider the remaining subassignment $h_0$ defined by 
$3\ord(t)=\ord(x)$.  
It breaks up into two subassignments, which we will call  
$h_c$ and $h_{nc}$ (for ``cubes'' and ``non-cubes'', respectively) 
defined, respectively, by $\exists y:y^3=x$ and 
$\nexists y: y^3=x$.
We omit $h_{nc}$, because it is similar to $h_1$, and focus on the most 
interesting part $h_c$.

We will use three extra residue field variables: the variable $\eta_1$ 
will stand for  for $\ac(x)$, $\eta_2$ for 
$\ac(t)$, and $\eta_3$ -- for the angular component of the difference between
$t$ and a given cube root of $x$ (the details will appear below, see equation 
(\ref{eq:pres1})).  We will also have one value sort variable $\gamma$ -- 
for the order of congruence between 
$t$ and the chosen cube root of $x$.
  
Now let us do this formally. We can take as the basis  the subassignment $C_1$ 
of $h[1,3,1]$ defined by the formula
\begin{multline}\label{eq:C_1}
%\begin{aligned}
\phi(x,\eta_1,\eta_2, \eta_3, \gamma)=\\
\text{`}(\exists y: y^3=x) 
\wedge  (\eta_1=\ac(x))\wedge (\eta_2^3=\ac(x))
\wedge (\gamma\ge \ord(x)+1)
\text{'}.
%\end{aligned}
\end{multline}
Let the function $c_1:C_1\to \A_{k((t))}^1$ be defined by 
$c_1(x,\eta_1,\eta_2,\eta_3, \gamma)=y$, where $y^3=x$ and $\ac(y)=\eta_2$. 
Note that $c_1$ is a definable function, since its graph clearly is a 
definable set.
Let the function $\alpha_1:C_1\to \Zgordon$ be defined by 
$\alpha_1(x,\eta_1,\eta_2, \eta_3,\gamma)=\gamma$, and let 
$\xi_1(x, \eta_1,\eta_2, \eta_3,\gamma)=\eta_3$.
We make the set $Z_{C_1,c_1,\alpha_1,\xi_1}$ with these data according to 
Definition \ref{def:cells}.
The presentation $\lambda:h_c\to Z_{C_1,c_1,\alpha_1,\xi_1}$ 
is given by the formula
\begin{equation}\label{eq:pres1}
%\begin{aligned}
\lambda(x,t):=(x,\ac(x),\ac(t), \ac(t-y),\ord(t-y)),
%\end{aligned}
\end{equation}
where $y^3=x$ and $\ac(y)=\ac(t)$.  Finally, let $\psi_1:C_1\to A$ be the function 
$$\psi(x,\eta_1,\eta_2, \eta_3, \gamma)=\lef^{-2\ord(x)-\gamma}. 
$$
\end{example} 
It is easy to see that all the conditions of cell decomposition theorem are satisfied with these formal constructions. We will soon see how this prepares the ground for integration, and will help us recover the calculation of Section~\ref{subsub:example}.

%%%%%%%%%%%%%%%%%%%%%%%%%%%%%%%%%%%%%%%%%%%%%%%%%%%%%%%%%%%
\subsection{Motivic integration as pushforward} We are almost ready to 
define integration with respect to one valued field variable. We just need 
to discuss the (tautological) 
integration with respect to extra residue field variables, 
and summation over $\Zgordon$-variables, since as we have just seen, 
we do pick up these variables in the process of 
cell decomposition.

\subsubsection{Integration over the residue field variables}\label{subsub:ras.variables}
Everything in this subsection comes from \cite[Section 5.6]{CL}. 

Let $f:S[0,n,0]\to S$ be the projection onto the first factor.
Recall that by definition, the ring of constructible functions on 
$S[0,n,0]$ is spanned by 
the elements of the form $a\otimes\varphi$, where $a$ is an element of 
$K_0(\rde_S[0,n,0])$, and $\varphi$ is a  
function on $S[0,n,0]$ with values in the ring $A$.
Using quantifier elimination, one can prove \cite[Proposition 5.3.1]{CL}  
that in fact it is enough to have just the elementary tensors of the form
$a\otimes \varphi$  where the  
$\varphi$'s are pullbacks to $S[0,n,0]$ of functions on $S$, namely,
the natural map 
\begin{equation}\label{eq:fns}
K_0(\rde_{S[0,n,0]})
\otimes_{{\mathcal P}^0(S)} {\mathcal P}(S) \to {\mathcal C}(S[0,n,0])
\end{equation}
is an isomorphism.

Here is an example illustrating this fact.
\begin{example}
Let $\varphi={\bf 1}_Y$ be a characteristic function of a definable 
subassignment $Y$ of $S[0,n,0]$. Then $Y$ is an element of $\rde_S$, so 
clearly ${\bf 1}_Y$ is in the image of the map (\ref{eq:fns}).
\end{example}

Given this isomorphism of rings of constructible motivic functions,
pushforward for the projection  $f$ is easy to define, and it is, essentially, tautological. An element $a$ of 
$K_0(\rde_S[0,n,0])$ can be viewed as an elements of $K_0(\rde_S)$ via 
composition of the map to $S[0,n,0]$ with $f$. We denote it by 
$f_{!}(a)$. Then let 
$f_{!}(a\otimes\varphi):=f_{!}(a)\otimes \varphi$. 

\subsubsection{Integration over $\Zgordon$-variables,\cite[Section 4.5]{CL}}\label{subsub:z.variables}
Essentially, the measure on $\Zgordon^r$ is just the counting measure, and integration is summation. More precisely, we call a  family $(a_i)$ of elements of 
$A$ {\bf summable}, if $\sum_i\nu_q(a_i)$ converges for all $q>1$. 
A function $\varphi(s,i)\in {\mathcal P}(S\times \Zgordon^r)$ is called 
{\bf $S$-integrable} if, for every $s\in S$, the family $(\varphi(s,i))_{i\in \Zgordon^r}$
is summable (recall that our functions are $A$-valued).
\begin{thm}\cite[Theorem-Definition 4.5.1]{CL}
For each $S$-integrable function $\varphi$ on $S\times \Zgordon^r$, there 
exists a unique constructible motivic function $\mu_S(\varphi)$ on $S$
such that for all $q>1$ and all $s$ in $S$,
$$
\nu_q(\mu_S(\varphi)(s))=\sum_{i\in \Zgordon^r}\nu_q(\varphi(s,i)).
$$ 
\end{thm}

The proof of this theorem requires cell decomposition for Presburger functions; we will not discuss it here.
One of the consequences of the structure of Presburger functions is the fact 
that the ring $A$ is the correct ring of values for constructible motivic functions. More precisely, it is the structure of Presburger functions that is ultimately responsible for the fact that it is enough to invert $\lef$ and the 
elements $1-\lef^{-n}$ in order to do integration of summable functions.

\subsubsection{Integration over a $1$-cell}\label{subsub:val.variables}
Let $S$ be a definable subassignment as before, and let 
$\pi:S[1,0,0]\to S$ be the projection onto the first factor.
Let $\varphi$ be a constructible motivic function on $S[1,0,0]$.
We want to produce a constructible motivic function $\pi_{!}(\varphi)$ 
on $S$ 
that is the result of integrating $\varphi$ along the fibers of $\pi$.
The idea of integration is very simple: take a cell decomposition 
of $S[1,0,0]$ adapted to $\varphi$. We have $S[1,0,0]=\sqcup Z_i$, 
where $Z_i$ are cells. 
The function $\varphi$ breaks up into the sum of its restrictions to cells:
$\varphi=\sum\varphi{\bf 1}_{Z_i}$, and we define the function 
$\pi_!(\varphi)$ cell by cell. 
If we care only for functions defined almost everywhere, we can discard the restriction of $\varphi$ to the union of $0$-cells, since it is supported on the set of smaller dimension than the restriction of $\varphi$ to the union of $1$-cells. 

Now let us define the pushforward on $1$-cells.
Let $Z$ be a $1$-cell, and let $\varphi_Z$ be the restriction of 
$\varphi$ to $Z$.
We have:
\[
\xymatrix{
& \ar[dd]^{p_1} Z_{C_1,\alpha,\xi,c} \ar[rrr]^{\hookrightarrow} &&& S[1,s,r]\ar[d]\\
& & \ar[ul]^{\lambda} Z \ar[rr]^{\hookrightarrow} && S[1,0,0]\ar[d]^{\pi} \\
A & \ar[l]^{\psi_1} C_1 \ar[r]_{j_1}^{\hookrightarrow} & S[0,s,r] \ar[r]_{\pi_1} & S[0,0,r] \ar[r] & S
}
\]

Note that by definition of the cell, $\varphi_Z$ is constant 
on the fibres of $p_1\circ\lambda$. If we identify $Z$ with 
$Z_{C_1,\alpha,\xi,c}$  by means of the presentation $\lambda$, we can pretend 
that the function $\varphi_Z$ lives on $Z_{C_1,\alpha,\xi,c}$, and it is 
constant  on the fibres of the projection  $p_1:Z_{C_1,\alpha,\xi,c}\to C_1$. 
It is natural to define the volume of the fibre of the projection $p_1$
over a point $y\in C_1$ to be  $\lef^{-\alpha_1(y)-1}$ -- by analogy with
the $p$-adic situation.
Hence, the following definition is natural:
\begin{definition}
\begin{equation}\label{eq:integral}
\pi_!(\varphi_Z):=\mu_S({\pi_1}_!({j_1}_!(\lef^{-\alpha_1-1} \psi_1))).
\end{equation}
\end{definition}   
Note that this definition automatically introduces normalization of the 
measure: by specifying the factor $\lef^{-\alpha_1-1}$, we have 
fixed the volumes of $1$-dimensional $p$-adic balls.

\begin{example}
Let us return to the example of section \ref{subsub:example}, and consider 
the case when the parameter $x$ is a cube, and in this case, let us only do the  integral over  a subset of the set $\{t\mid 3\ord(t)=\ord(x)\}$. 
Recall the notation from section \ref{subsub:example}: we had the set $B_1$ of all $t$ that are close to the cube root (or one of the three cube roots) of $x$. In Example~\ref{ex:cells}, we defined the corresponding subassignment 
$h_c$ and 
showed that it is a $1$-cell in the cell decomposition adapted to the 
constructible function $\varphi(x,t)=\lef^{-\ord(t^3-x)}$. 
Let us now compute the motivic integral of $\varphi$ with respect to the 
variable $t$ over the cell $Z=h_c$. 

In the notation used in the above definition,
we have $S=\A^1_{k((t))}=h[1,0,0]$.   
On $Z=h_c$ (see Example~\ref{ex:cells}), we have 
$\varphi=\lambda^{\ast}p_1^{\ast}(\psi_1)$, where $\psi_1$ is a function 
on the basis $C_1\subset h[1,3,1]$ defined by 
$\psi_1(x, \eta_1,\eta_2,\eta_3,\gamma)=\lef^{-2\ord(x)-\gamma}$.
The function $\alpha_1$ on $C_1$ is defined by:
$\alpha_1(x,\eta_1,\eta_2,\eta_3,\gamma_2)=\gamma$, so we have
$\lef^{-\alpha_1-1} \psi_1=\lef^{-2\ord(x)-2\gamma-1}$.
Note that this function is in ${\mathcal P}(C_1)$. 
By definition, 
${\pi_1}_{!}{j_1}_{!}(\lef^{-\alpha_1-1} \psi_1)=
[C_1]\otimes\lef^{-2\ord(x)-2\gamma-1}$, where now $C_1$ is thought of 
as an element of $\rde_{S[0,0,1]}$ via the map $\pi_1\circ j_1$, and 
$[C_1]$ is its class in $K_0(\rde_{S[0,0,1]})$.
Let us denote the projection (that forgets the $\Zgordon$-variable)
from $S[0,0,1]$ to $S$ by $p$. 
Now, $\mu_S$ amounts to summation over $\gamma$, and we get
\begin{equation}\label{eq:pi!}
\pi_!(\varphi_Z):=\mu_S({\pi_1}_!({j_1}_!(\lef^{-\alpha_1-1} \psi_1)))=
[p(C_1)]\otimes\lef^{-2\ord(x)-1}\lef^{-2(\ord(x)+1)}(1-\lef^{-2}).
\end{equation}

Recall that $C_1$ is defined by the formula (\ref{eq:C_1}) of 
Example~\ref{ex:cells}. Then $p(C_1)$ is a subassignment of $S[1,3,0]$
defined by: $\eta_1=\ac(x), \eta_2^3=\eta_1$ (we call the three residue field 
variables $\eta_{1,2,3}$). 
Note that magic happens as we fix a local field $K$ with a uniformizer 
$\varpi_K$ and residue field $\F_q$, and interpret all the formulas in it.
As we discussed briefly in Section~\ref{sub:constr.f} and as we shall see in 
detail in Section~\ref{sec:back}, to make $[p(C_1)]$ into a function on $S$,
we just need to count, for $x\in S(K)$, the number of points on the fibre of
$C_1$ over $x$. In our case, this yields three possible values of $\eta_2$ for each fixed $\eta_1=\ac(x)$ if there are $3$ cube roots of $1$ in the field, 
or just one value of $\eta_2$ if there is only one cube root.
Since $\eta_3$ can take any value except $0$, we get $3(q-1)$ or $q-1$, 
respectively. If we plug these numbers into (\ref{eq:pi!}), and replace all occurrences of $\lef$ with $q$, we get an answer that agrees with equation 
(\ref{eq:B_1}) of section~\ref{subsub:example}.
\end{example}

\subsection{What was swept under the carpet}\label{sub:carpet}
Since our goal was just to give a very basic exposition of the main ideas 
of the theory of motivic integration, we have left out, so far, some very 
important issues, such as integrability and integration over manifolds.
 
\subsubsection{Integrability}
Naturally, there are many definable sets whose $p$-adic volume is not finite, 
and there are many constructible motivic functions whose integral 
should not converge. In the earlier versions of motivic integration this 
issue was mainly dealt with by letting the valued field variables in all 
formulas range only  over the ring of integers, and not over the whole 
valued field. That approach made the domain of integration compact, and 
guaranteed finiteness of the volume.

One of the advantages of the theory developed in \cite{CL} is that the restriction to the ring of integers is dropped, and instead a natural class of 
integrable functions is constructed. This is done by 
starting out only with {\emph {summable}} 
Presburger functions over $\Zgordon^r$; 
as the valued-field and residue-field variables are added, it is necessary to 
consider Grothendieck {\emph {semirings}} of the so-called {\emph {positive}} 
constructible motivic functions, instead of the full rings of constructible motivic functions.
Essentially, the term ``positive'' comes from the partial order that we have on the ring of values $A$. 
The semiring of positive constructible motivic Functions on $S$ is 
denoted by $C_+(S)$.
 We refer to \cite[5.3]{CL} or 
\cite[3.2]{CL.expo} for details.
The class of integrable positive Functions on $Z\in \de_S$ 
(denoted by $I_S C_+(Z)$) is defined inductively along with 
the procedure of integration itself.

Let $S$ be in $\de_k$. 
The main existence theorem for motivic integral 
(Theorem \cite[Theorem 10.1.1]{CL}) states that 
there is a unique functor from the category 
$\de_S$ to the category of abelian semigroups, $Z\mapsto I_S C_{+}(Z)$, assigning to every morphism $f:Z\to Y$ in $\de_S$ a morphism 
$f_{!}:I_SC_+(Z)\to I_SC_+(Y)$ that satisfies a list of axioms.
We have already discussed most of these axioms in some form: they 
include additivity, natural behaviour with respect to inclusions and 
projections, normalization according to (\ref{eq:integral}), and the
Jacobian transformation rule, which is discussed below in \ref{subsub:graphs}. 
Note that pushforward is functorial,
in the sense that it respects compositions: $(f\circ g)_{!}=f_{!}\circ g_{!}$.
We refer to \cite[Theorem 10.1.1]{CL} or to \cite[Section 2.5]{CHL} for 
the complete list. 

\subsubsection{Integration over graphs}\label{subsub:graphs}
The idea of integration that we have sketched so far is sufficient for 
integration of constructible functions over $d$-dimensional subsets of 
$\A^d_{k((t))}$ for some $d$.  It would be natural for the theory to 
include integration over manifolds, and a Jacobian transformation rule. 
Cell decomposition helps with this issue as well: $0$-cells are 
basically graphs of functions, and so one can make sure that transformation 
rule holds by defining integrals over $0$-cells appropriately. 

For a definable subassignment $h$, let 
${\mathcal A}(h)$ be the ring of definable functions from $h$ to 
$\A_{k((t))}^1$. For every positive integer $i$, one can define 
an ${\mathcal A}(h)$-module $\Omega^i(h)$ of definable $i$-forms on $h$.  
As one would naturally hope, the module of top degree forms is free 
of rank $1$, and there is a canonical volume form  
$|\omega_0(h)|$, which is an analogue 
of the canonical volume form in the $p$-adic case.

\begin{definition}\cite[8.4]{CL}
Let $f:X\to Y$ be a morphism between two definable subassignments 
of $h[m,n,r]$ and $h[m',n',r']$, respectively. Assume that both $X$ 
and $Y$ are of dimension $d$, and the fibres of $f$ have dimension $0$. 
Then the order of Jacobian
\footnote{In geometric motivic integration, the {\bf order of Jacobian}
is given a very geometric meaning: if $f:X\to Y$ is a morphism of varieties,
the order of Jacobian is the function on the arc space of $X$ that assigns 
to each arc its  order of tangency to the singular locus of the morphism $f$. 

As we discuss in Appendix 1, motivic integration theory described here 
specializes to 
geometric motivic integration. It is worth pointing out that the two 
notions of the order of Jacobian 
agree, \cite[8.6]{CL}.
}
 is defined naturally by the formula
\footnote{It is possible to show that a definable function on a definable subassignment  $S$ is analytic outside a subassignment $S'$ with $\dim S'< \dim S$.
On the subassignment $S\setminus S'$ the usual determinant formula for Jacobian holds. 
}
$$f^{\ast}|\omega_0|_Y=\lef^{-\text{ordjac} f}|\omega_0|_X,$$
with $\text{ordjac} f$ a $\Zgordon$-valued function on $X$ defined outside a
definable subassignment of dimension less than $d$.  
\end{definition}

Now let $Z$ be a $0$-cell that is part of cell decomposition adapted to a 
constructible fucntion $\varphi$, and let $\varphi_Z:=\varphi{\bf 1}_Z$ 
be the restriction of $\varphi$ to $Z$. Let us  assume here that $Z$ has 
dimension $d$, and this is the dimension of the support of $\varphi$.  
Then we have:
\[
\xymatrix{
& Z_{C_0,c} \ar[dd]^{p_0} \ar[rr]^{\hookrightarrow} && S[1,s,0] \ar[d] \\
&& \ar[ul]^{\lambda} Z \ar[r]^{\hookrightarrow} & S[1,0,0]\ar[d]^{\pi} \\
A &\ar[l]^{\psi_0} C_0 \ar[r]_{j_0}^{\hookrightarrow} & S[0,s,0] \ar[r]_{\pi_0} & S
}
\]

Recall that by definition of $\psi_0$, we have 
$\varphi_Z=\lambda^{\ast}p_0^{\ast}(\psi_0)$. As in the case of $1$-cells, 
let us imagine that $Z$ is identified with $Z_{C_0,c}$ by means of the 
isomorphism
$\lambda$, and the function $\varphi_Z$ is a function on $Z_{C_0,c}$. 
\footnote{Of course, when the construction is finished, one needs to prove that it does not depend on $\lambda$. This turns out to be the case, see \cite[\S~9.1-9.2]{CL}.}
By definition of a $0$-cell, the fibres of the projection $p_0$ are 
$0$-dimensional, so what we expect is that the functions $\varphi_Z$ and
$\psi_0$ would be related essentially by a factor that captures the order of the Jacobian 
of the map between $Z$ and $C$. This is exactly the case.  
By definition, $Z_{C_0,c}$ is an image of $C$ under the map $p_0^{-1}$.
It is natural to define ${p_0}_!(\varphi_Z)$ as $\lef^{\gamma}\psi_0$, where
the function $\gamma:C_0\to \Zgordon$ is defined by $y\mapsto (\text{ordjac}p_0)\circ p_0^{-1}$.
Finally, we already know how to define ${\pi_0}_{!}$ 
(subsection \ref{subsub:ras.variables}), 
and ${j_0}_!$ (extension by zero).
Putting all these pieces together, we get
\begin{definition}
$$
\pi_!(\varphi_Z):={\pi_0}_!({j_0}_!(\psi_0\lef^{\gamma})).
$$
\end{definition}

The harderst part of the theory is proving that the final definition 
of pushforward does not depend on the choice of cell decomposition, and that 
integration with respect to several 
valued field variables does not depend on the order 
(statements of Fubini type).

\subsection{Motivic volume}
Let $\Lambda\in \de_k$ be a definable subassignment, and let 
$S\in \de_{\Lambda}$ (in particular, $S$ comes equipped with a morphism
$f:S\to\Lambda$). Then we can define {the} relative motivic volume of $S$ as
$$\mu_{\Lambda}(S)=f_{!}([{\bf 1}_S]).$$ 
In particular, when $\Lambda=h_{\spec k}$ is the final object of the 
category $\de_k$, we get the motivic volume 
for all definable subassignments $S$ such that the characteristic 
function ${\bf 1}_S$ is integrable.

Let us call a subassignment $Z$ of some $h[m,n,0]$ bounded if there exists
a positive integer $s$ such that $Z$ is contained in the subassignment $W_s$ of
$h[m,n,0]$ defined by $\ord(x_i)\ge -s$, $1\le i\le m$ 
(where the variables $x_i$ run over the  valued field). 

\begin{proposition}\cite[Proposition 12.2.2]{CL}
If $Z$ is a bounded definable subassignment of $h[m,n,0]$, then 
$[{\bf 1}_Z]$ is integrable.
\end{proposition}

By definition, motivic volume takes values in the ring of positive integrable 
constructible 
motivic Functions on $\spec k$. This ring, by definition,  
is 
$$
SK_0(\rde_{\spec k})\otimes_{\N[\lef-1]}A_{+},
$$
where $SK_0(\rde_{\spec k})$ is the Grothendieck semiring (as opposed to 
the full Grothendieck ring) that is made by taking only formal 
linear combinations of equivalence classes of objects of $\rde_{\spec k}$ with
nonnegative coefficients, and $A_+$ is the set of nonnegative elements of 
$A$.

\section{Back to $p$-adic integration}\label{sec:back}

Everywhere in this section, we fix the base field $k=\Qgordon$, for simplicity of the exposition.
Recall that in the definition of the language of Denef-Pas (see Section 
\ref{sub:ldp}), there was some flexibility in the matter of choosing what to add to the language as allowed coefficients for formulas.
\emph{Everywhere in this section, we will consider one specific variant of Denef-Pas language: we allow coefficients in $\Zgordon[[t]]$ for the valued field sort, and coefficients in $\Zgordon$ for the residue field sort. This language will be denoted
${\mathcal L}_{\Zgordon}$.}

There are two collections of fields {over which} we would like to do integration:
local fields of characteristic zero, and the function fields $\F_q((t))$.
Let ${\mathcal A}_{\Zgordon}$ be the collection of all finite field 
extensions of non-archimedean completions of 
$\Qgordon$, and let ${\mathcal B}_{\Zgordon}$ be the collection of all local fields of positive characteristic. 

In the last section we sketched the construction of a motivic volume of a 
subassignment $h\in \de_{\Qgordon}$, and more generally, of an integral of a constructible motivic function on $h$. In order to relate this motivic integration 
with the classical $p$-adic integration of Section~\ref{sub:Weil}, we need 
to do two things: first, we need to relate subassignments to the $p$-adic measurable sets, and second, we need to find a way to get from the values of the motivic volume to the rational numbers. 
We start with the first task.

\subsection{Interpreting formulas in $p$-adic fields}\label{sub:interpr}

Observe that a definable subassignment $S$ of, say, $h[1,1,0]$ does not automatically give us a subset of $\Qgordon_p\times \F_p$: indeed, $S(\Qgordon_p)$ is by definition a 
subset of 
$\Qgordon_p((t))\times \Qgordon_p$ rather than of $\Qgordon_p\times \F_p$.
However, is is clear that we can interpret the formulas defining $S$ so that we would get a subset of $\Qgordon_p\times \F_p$ as desired.
Let us describe this procedure precisely (we are essentially quoting 
\cite[\S~6.7]{CL.expo}).

Let $S$ be a definable subassignment of $h[m,n,r]$. As specified at the beginning of this section, by this we mean that $S$ is defined by a formula $\varphi$ 
in 
Denef-Pas language with coefficients in $\Zgordon[[t]]$. 
Let $(K,\varpi_K)$ be a local field of characteristic $0$ from the collection 
${\mathcal A}_{\Zgordon}$, with the choice of a uniformizer.  
The field $K$ can be considered as a $\Zgordon[[t]]$-algebra via the morphism
$$\lambda_{\Zgordon, K}:\Zgordon[[t]]\to K: \quad \sum_{i\ge 0}a_i t^i\mapsto 
\sum_{i\ge 0} a_i\varpi_K^i.$$
Note that the series $\sum_{i\ge 0} a_i\varpi_K^i$ converges in $K$, since 
$\ord(a_i)\ge 0$ for any $a_i\in \Zgordon$.
A similar morphism exists also for fields of finite characteristic from 
the collection ${\mathcal B}_{\Zgordon}$, even though in this case we prefer to 
write it as
$$\lambda_{\Zgordon, K}:\Zgordon[[t]]\to K: \quad \sum_{i\ge 0}a_i t^i\mapsto 
\sum_{i\ge 0} (a_i\pmod{p_K})\varpi_K^i,$$
where $p_K$ is the characteristic of the residue field of $K$.

Using these  morphisms, any formula $\varphi$ with coefficients in $\Zgordon[[t]]$ and $m$ free variables of the valued field sort and no other free variables
can be interpreted to define a subset of $\A^m(K)$ for any 
$K\in {\mathcal A}_{\Zgordon}\cup {\mathcal B}_{\Zgordon}$. 
Formulas in the language of rings with coefficients in $\Zgordon$ can naturally be 
interpreted in the residue field of $K$, via reduction $\mod q_K$.
There is no additional work needed for the variables running over $\Zgordon$.
This way, any definable (with the mentioned above restriction on coefficients)
subassignment $S$ of $h[m,n,r]$ gives a subset $S_{K,\phi}$ of 
$K\times k_{K}\times \Zgordon^r$, where 
$K\in {\mathcal A}_{\Zgordon}\cup{\mathcal B}_{\Zgordon}$, and $k_K$ is the residue field of $K$, and where $\phi$ is the formula (or collection of formulas) 
defining the subassignment $S$.
  
There is a very important issue here: the set $S_{K, \phi}$
depends on the choice of the formula $\phi$ that we used to define $S$, 
as illustrated by a very simple example.
Consider the two formulas 
$\phi_1(x)=\text{`}{x=0}\text {'}$ and $\phi_2(x)=\text{`}{3x=0}\text {'}$. 
For each field $K$ of characteristic $0$, either formula defines a one-point 
set $\{0\}$, so $\phi_1$ and $\phi_2$ define the same subassignment 
(call is $S$) of 
$h[1,0,0]$. On the other hand, for the fields $K$ of characteristic $3$, 
$S_{K,\phi_1}\neq S_{K,\phi_2}$.
This example illustrates that 
\emph{the correspondence between definable subassignments and definable $p$-adic sets is well defined only for sufficiently large $p$. Moreover, the choice of the primes to discard depends on the formula we are using to describe a given set, not on the set itself.}
The fact that only finitely many primes need to be discarded (which is of course crucial) is a nontrivial theorem.
Precisely, we have:

\begin{proposition}\cite[\S\S~6.7, 7.2]{CL.expo}
If two formulas $\psi$ and $\psi'$ define the same subassignment $S$, 
then there exists an integer $N$ such that 
$S_{K,\psi}=S_{K,\psi'}$ for every field 
$K\in {\mathcal A}_{\Zgordon}\cup {\mathcal B}_{\Zgordon}$
with residue characteristic greater or equal to $N$.
However, this number $N$ can be arbitrarily large for different $\psi'$.
\end{proposition}

\subsubsection{Specialization of constructible motivic Functions}\label{sub:specialization}
We have just described how definable subassignments give measurable subsets of 
$p$-adic fields. Let us now describe the specialization of constructible 
motivic functions.

First, note that a morphism of definable subassignments $f:Z\to W$ specializes
to a function $f_K:S_K\to W_K$ for all 
$K\in {\mathcal A}_{\Qgordon}\cup{\mathcal B}_{\Qgordon}$ of sufficiently large 
residue 
characteristic (since the graph of $f$ is a definable subassignment, it 
will specialize to a definable subset of $S_K\times W_K$, and that gives the graph of $f_K$).
In particular, for $S\in \de_{\Qgordon}$, $\Zgordon$-valued functions on $S$ specializes
to $\Zgordon$-valued functions on $S_K$. The functions with values in the ring 
$A$ specialize to $\Qgordon$-valued functions once we replace $\lef$ with $q$, 
where $q$ is the cardinality of the residue field of $K$. 
Thus we can interpret elements of ${\mathcal P}(S)$.
 
Recall that a constructible motivic function on $S$ is an element of 
${\mathcal P}(S)\otimes K_0(\rde_S)$.
As mentioned in Section~\ref{sub:constr.f}, an 
element $\pi:W\to S$ of $K_0(\rde_S)$ gives an  integer-valued functions 
on $S_K$ by $x\mapsto \#\pi_K^{-1}(x)$.  

The main point is that motivic integration specializes to $p$-adic 
integration. Since now we also have the residue-field and integer-valued parameters, when we consider $p$-adic measure, we take the product of 
Serre-Oesterl\'e measure on the Zariski closure of the set cut out by 
the valued-field variables with the counting measure on 
$k_K^n\times \Zgordon^r$.

Let $\Lambda\in \de_{\Qgordon}$ be a definable subassignment.
Let $S\in \de_{\Lambda}$, with the morphism $f:S\to \Lambda$.
Let $\varphi$ be an integrable constructible motivic function on $S$, and let $K$ be a local field. Then we have $f_K:S_K\to \Lambda_K$, and the 
interpretation $\varphi_K$, which is a function on $S_K$ (all these
are well defined when the residue characteristic of $K$ is large enough).
It is possible to prove that the restriction of $\varphi_K$ to the fibre of 
$f_K$ at a point $\lambda\in \Lambda_K$ is integrable for almost all 
$\lambda\in \Lambda_K$. We denote by $\mu_{\Lambda_K}(\varphi_K)$ the 
function on $\Lambda_K$ that assigns to each point $\lambda\in \Lambda_K$ 
the integral of 
$\varphi_K$ over the fibre of $f_K$ at $\lambda$.

\begin{theorem}\cite[9.1.5, Specialization Principle]{CLF}\label{thm:spec.pr}
Let $f:S\to \Lambda$ be an ${\mathcal L}_{\Zgordon}$-definable morphism, and let
$\varphi$ be a constructible motivic function on $S$, relatively integrable
with respect to $f$.  
Then there exists $N>0$ such that for all $K$ in 
${\mathcal A}_{\Zgordon}\cup{\mathcal B}_{\Zgordon}$ with residue characteristic 
greater than $N$, and
every choice of the uniformizer $\varpi$ of the valuation on $K$, 
$$
(\mu_{\Lambda}(\varphi))_K=\mu_{\Lambda_K}(\varphi_K).
$$           
\end{theorem} 

This theorem is proved by comparing the construction of the motivic integral 
with the understanding of the $p$-adic measure that one gets from $p$-adic 
cell decomposition theorem \cite{Denef}. 

\subsection{Pseudofinite fields}\label{sub:pseudo}

By now we have the motivic volume with values in 
$SK_0(\rde_{\spec k})\otimes_{\N[\lef-1]} A_+$, and it specializes to the classical $p$-adic volume for almost all $p$, as discussed above.
It turns out that if we just want to capture the $p$-adic volume, then our 
motivic volume is a bit too refined and complicated object, namely, one can 
identify a lot of elements of $K_0(\rde_{\spec k})$, and specialization would 
still hold.  In order to define a new equivalence relation on formulas in the 
language of rings, we need to define the category of 
pseudofinite fields first.

\begin{definition}
The field $K$ of characteristic zero is called {\bf pseudofinite} if it is perfect, has exactly one field extension of
each finite degree, and if $V$ is a geometrically irreducible variety over $K$, then $V$ has a 
$K$-rational point.
\end{definition}

One can get an example of a pseudofinite field by means of 
constructing an {\emph {ultraproduct}} of finite fields, see \eg, 
\cite[\S~20.10]{FJ}.

\begin{definition}\cite{DL.congr}.
Let {\bf $K_0(\pff_k)$} be the group generated by symbols $[\phi]$, where
$\phi$ is any formula in the language of rings over $k$, subject to the 
relations:
$[\phi_1 \vee \phi_2]=[\phi_1]+[\phi_2]-[\phi_1\wedge \phi_2]$ whenever 
$\phi_1$  and $\phi_2$ have the same free variables, and the relations 
$[\phi_1]=[\phi_2]$
if there exists a ring formula $\psi$  over $k$ such that the interpretation of $\psi$ in 
any pseudofinite field $K$ containing $k$ gives a graph of a bijection between the tuples of elements of  $K$ satisfying $\phi_1$ and those satisfying 
$\phi_2$.   
The multiplication on $K_0(\pff_k)$ is induced by the conjunction of 
formulas in disjoint sets of variables.
The additive group of $K_0(\pff_k)$  is called the Grothendieck group of 
pseudofinite fields.
\end{definition}

The reason the category of  pseudofinite fields turns out to be so useful 
for us is the following theorem. 
A DVR-formula is a formula in the language of Denef-Pas with coefficients in 
$\Zgordon[[t]]$ in the {valued} field sort, and such that all its valued field 
variables are restricted to the ring of integers (DVR stands for ``Discrete Valuation Rings''). 

\begin{theorem} (Ax-Kochen-Ersov Principle)
Let $\sigma$ be a DVR-formula over $\Zgordon$ with no free variables. Then the following statements are equivalent:
\begin{enumerate}
\item The interpretation of $\sigma$ in $\Zgordon_p$ is true for all but finitely many primes.
\item The interpretation of $\sigma$ in $K[[t]]$ is true for all pseudofinite fields $K$.
\end{enumerate}
\end{theorem}
 
\subsection{Comparison theorems}\label{sub:comp}
Let the base field be $k=\Qgordon$, as before.
Given a definable subassignment $X$ of $h[m,0,0]$, by now we have defined 
the associated with it subsets of $K^m$ for 
$K\in {\mathcal A}_{\Zgordon}\cup {\mathcal B}_{\Zgordon}$, and we have defined the 
motivic  volume of $X$,  $\mu(X)\in SK_0(\rde_{\Qgordon})\otimes_{\N[\lef-1]}A_+$.
There is a natural map from the Grothendieck ring  $K_0(\rde_k)$ 
to $K_0(\pff_k)$: we just
identify the subassignments that coincide on the category of pseudofinite 
fields containing $k$ to obtain a class in $K_0(\pff_k)$.
Hence, to each subassignment $X$ we have also associated an element 
of $K_0(\pff_{\Qgordon})$, which we will also denote by $\mu(X)$. 

By Ax-Kochen-Ersov principle, two formulas
$\phi_1$ and $\phi_2$ define the subsets of $K^m$ of the same volume  for 
$K\in {\mathcal A}_{\Zgordon}\cup {\mathcal B}_{\Zgordon}$ with residue characteristic 
bigger than $N$ for some $N$ if and only if 
$\mu(h_{\phi_1})=\mu(h_{\phi_2})$, where $h_{\phi}$ denotes the 
subassignment defined by the formula $\phi$.  

It is a difficult theorem of Denef and Loeser \cite{DL} that  
there exists a unique ring morphism
$$\chi_c:K_0(\pff_k)\to K_0^{mot}(\var_k)\otimes \Qgordon,$$
that satisfies two natural conditions. The first condition is that for any formula $\varphi$ which is a conjunction of polynomial equations over $k$, the element $\chi_c([\varphi])$ equals the class in $K_0^{mot}(\var_k)\otimes \Qgordon$
of the variety defined by $\varphi$. The seconds condition is more complicated: it specifies how the map $\chi_c$ should behave with respect to cyclic covers. This relates to elimination of quantifiers in formulas of the form  
$\varphi(x)=\text{'}\exists y: y^d=x\text{'}$. 
It is this condition that makes Chow motives the right category for the values of the volume, as opposed to varieties, which would not have been 
sufficient. 

We refer to \cite[Th.~2.1]{DL.congr}
for the precise statement and a sketch of the proof, and to \cite{Tom.intro}
for an exposition.

The existence of the map $\chi_c$ allows to state the Comparison Theorem, 
\cite[Th.~8.3.1, Th.~8.3.2]{DL.arithm}.
Here we quote a reformulation of this theorem as stated in \cite{CL.expo}.
\begin{theorem}
Let $\varphi$ be a formula in the language of Denef-Pas, with $m$ free 
valued field variables and no other free variables. There exists a 
virtual motive $M_{\varphi}$, canonically attached to $\varphi$, such that, 
for almost all prime numbers $p$, the volume of $h_{\Qgordon_p,\varphi}$ is finite 
if and only if the volume of $h_{\F_p[[t]],\varphi}$ is finite, and in 
this case they are both equal to the number of points 
of $M_{\varphi}$ in $\F_p$.\footnote{
In the original construction, the virtual Chow motive $M_{\varphi}$ lives in a certain completion of the ring $K_0^{mot}(\var_k)$ (see {Appendix}~1, and 
\cite{Tom.intro}). It follows from the Cluckers-Loeser theory 
of motivic integration described in the previous section 
that $M_{\varphi}$ lives in the ring obtained 
from $K_0^{mot}(\var_k)\otimes \Qgordon$ by inverting $\lef$ and $1-\lef^{-n}$ for 
all positive integers $n$.
When we say ``the number of {points} on $M_{\varphi}$'' we mean by this the 
extension of the function that counts the number of points over $\F_q$ from the category of varieties to the ring where $M_{\varphi}$ lives. This 
extension is {obtained} as follows: first, one replaces the number of points by the 
alternating sum of the trace of Frobenius on cohomology 
(as in Grothendieck-Lefschetz fixed point formula). This procedure is well-defined for Chow motives, and extends the notion of the number of rational 
points of a variety. Then the Trace of Frobenius function is extended 
to the Grothendieck ring by additivity, 
and then extended further to the tensor product with $\Qgordon$, in a 
natural way (at this point it becomes $\Qgordon$-valued).
Finally, if we assign the value $q$ to $\lef$, this function extends to the 
localization by $\lef$ and $1-\lef^{-n}$.
}
\end{theorem}

\begin{remark}
Even though it is necessary to make a map from $K_0(\rde_k)$ 
to $K_0(\text{PFF}_k)$ and further
to the  ring of virtual Chow motives in order to state the comparison theorems
that give a  geometric interpretation of the $p$-adic measure, 
the motivic volume taking values in $SK_0(\rde_{\spec k})\otimes A_+$ is sufficient for 
the transfer principle that we state in the next section. In fact, Ax-Kochen-Ersov principle that we have referred to in order to justify the map to $K_0(\text{PFF}_k)$ follows from this general transfer principle.
The way to think about it is that the motivic volume in 
$SK_0(\rde_k)\otimes A_+$ is the finest invariant of a subassignment; depending on the context, 
one can map it to more crude invariants.
For example, motivic integration specializes to integration with respect to 
Euler characteristic, as explained in the Introduction to \cite{CL}; one can also get Hodge or Betti numbers from the motivic volume 
(that was one of the first applications of motivic integration), and so on.
\end{remark}

\section{Some applications}\label{sec:applic}

There are two natural directions for application of arithmetic motivic 
integration. One is, to get various ``uniformity in $p$'' results. 
A very spectacular application in this direction is 
the results of Denef and Loeser on rationality of Poincar\'e series.
There are excellent expositions 
\cite{DL.congr} and \cite{Denef}, so we will not discuss it here.

The other direction is transfer of identities from function fields to fields of characteristic zero. This is made possible by the very general transfer principle, which follows immediately from the construction of the 
motivic integral and the fact that it specializes to the $p$-adic integral.

\begin{theorem}\cite[Transfer principle for integrals with parameters.]{CLF}
Let $S\to \Lambda$ and $S'\to\Lambda$ be 
${\mathcal L}_{\Zgordon}$-definable morphisms.
Let $\varphi$ and $\varphi'$ be  ${\mathcal L}_{\Zgordon}$-constructible 
motivic functions on $S$ and $S'$, respectively. There exists $N$ such that
for every $K_1$ in ${\mathcal A}_{\Zgordon, N}$ and $K_2$ in 
${\mathcal B}_{\Zgordon,N}$ 
with isomorphic residue fields,
$$\mu_{\Lambda_{K_1}}(\varphi_{K_1})=\mu_{\Lambda_{K_1}}(\varphi_{K_1}')
\quad\text{if and only if}\quad
\mu_{\Lambda_{K_2}}(\varphi_{K_2})=\mu_{\Lambda_{K_2}}(\varphi_{K_2}').
$$
\end{theorem}
Loosely speaking, this theorem says that an equality of 
integrals of the specializations of 
two constructible motivic functions holds  over all local fields 
of characteristic zero with sufficiently large residue characteristic 
if and only if it holds over all function fields with sufficiently large 
residue characteristic. 

The most recent and important application of this transfer principle is 
the transfer principle for the Fundamental Lemma that appeared in \cite{CHL}.
Here we cannot {explain} the {Fundamental} Lemma (which states that certain 
$\kappa$-orbital integrals on two related groups are equal), 
so we only include a
brief discussion of the relevance of motivic integration to computing 
orbital integrals.

\subsection{Orbital integrals} 
Recall the definition:
\begin{definition}
Let  $G$ be a $p$-adic group and $\lie$ -- its Lie algebra, and let 
$X\in \lie$. 
An {\bf orbital integral at $X$} is a distribution on the space of Schwartz-Bruhat 
functions on $\lie$ defined by
$$\Phi_G(X,f):=\int_{G/C_G(X)} f(g^{-1}Xg)\id^{\ast}g,
$$
where $C_G(X)$ is the centralizer of $X$ in $G$, and $\id^{\ast}g$ is 
the invariant measure on $G/C_G(X)$.
\end{definition} 

The natural question (posed by T.C. Hales, \cite{Tom.p-adic}) is, 
can one use motivic integration to compute the orbital integrals in a 
$p$-independent way?  

Using all the terminology introduced above, we can rephrase this question:
\emph{Suppose we have fixed a definable test function $f$. 
Is the orbital integral $\phi_G(X,f)$ a constructible function of $X$?}
 
It looks like a constructible function, because we start with a definable 
function $f(g^{-1}Xg)$ of two variables $X$ and $g$, and then integrate 
with respect to one of the variables -- so by the main result of the 
theory of motivic integration, we should get a constructible function of the 
remaining variable. The difficulty, however, lies in the 
fact that the space of integration and the measure $d^{\ast}g$ on it vary 
with $X$.

The initial approach taken in \cite{Tom.orbital} and \cite{CH} was to 
average the orbital integral over definable sets of elements $X$, and then use 
local constancy results to make conclusions about 
the individual ones.

In \cite{CHL}, the authors start with definability of field extensions, 
(which leads to definability of centralizers),  and gradually prove that 
all ingredients of the definitions of the so-called $\kappa$-orbital 
integrals  appearing in the Fundamental Lemma are definable.
Consequently,
\begin{theorem}(Cluckers-Hales-Loeser, \cite{CHL})\label{thm:transfer}
The transfer principle applies to the Fundamental Lemma. 
\end{theorem}

It follows, in particular, from the main results of \cite{CHL} 
that the answer to our question is affirmative: $\Phi_G(X,f)$ is a 
constructible function of $X$ when $f$ is a 
fixed definable function.

We also observe that the results of \cite{CH} give quite precise information 
about the restriction of this  constructible function to the set of 
so-called {\bf good} elements. This direction is pursued further in 
\cite{CGS} with the hope of developing an actual algorithm for 
{computing} orbital integrals.

\subsection{Harish-Chandra characters}

Let $G$ be a $p$-adic group, and let $\pi$ be a representation of $G$. 
Harish-Chandra distribution character of $\pi$ is also defined as an 
integral over $G$, so it is natural to ask if the character is motivic as well.
The main difficulty in answering this question is that the construction of 
representations has many ingredients, and does not a priori appear to 
be a definable construction. However, if one adds additive characters of the 
field to the language (for example, by passing to the exponential functions 
as discussed in Section~\ref{sub:fourier}), then it is very likely that the 
construction of representations can be carried out within the language. 
Some partial results stating that certain classes 
of Harish-Chandra characters, when restricted to the neighbourhood of 
the identity, are motivic, appear in \cite{G} for depth-zero representations of classical groups, and in \cite{CGS} for certain positive depth 
representations.   

To give a flavour of a motivic calculation that appears when dealing with characters (and orbital integrals), 
we have included Appendix 2, where we compute 
motivic volume of a set that is relevant to the values of 
characters of depth zero representations of
$G=\text{SL}(2,K)$, where $K$ is a $p$-adic field. 
Many more calculations of this kind can be found in \cite{CG}.

\subsection{Motivic exponential Functions, and Fourier 
transform}\label{sub:fourier}
In \cite{CLF}, R. Cluckers and F. Loeser developed a complete theory of Fourier transform for the motivic measure described above.
Here we sketch the main features of this theory, since it is used in the proof 
of Theorem \ref{thm:transfer}, and is certain to find many 
other applications.

\subsubsection{Additive characters}\label{subsub:characters} 
We start by recalling the information about 
additive characters of valued and finite fields.

First, for a prime field $\F_p$, we can identify the elements of the field
with the integers $\{0,1,\dots,p-1\}$. Then one character of the additive group
of $\F_p$ can be written explicitly as 
$x\mapsto \exp\left(\frac{2\pi i}{p} x\right)$, and it generates the dual group  $\hat\F_p$ of $\F_p$.
For a general finite field $\F_q$ with $q=p^r$, we can explicitly  write down 
one character by composing our generator of $\hat\F_p$ with the trace map:
\begin{equation}\label{eq:generator}
\psi_0:x\mapsto \exp\left(\frac{2\pi i}{p}\text{Tr}_{\F_q/\F_p}(x)\right).
\end{equation}
This way to write the character allows us to talk about characters as ``exponential functions'', and this will be used in the next subsection. 
Given this character, we can identify the additive group of $\F_q$ with 
its Pontryagin dual via the map $a\mapsto\psi_0(ax)$.

The additive group of a local field $K$ is self-dual in a similar way.
If $\psi:K\to\Cgordon^{\ast}$ is a nontrivial character, then 
$a\mapsto \psi(ax)$ gives an isomorphism between $K$ and $\hat K$.

In particular, in agreement with our choice of the identification of $\F_p$ with $\{0,\dots, p-1\}$, and of $\psi_0$ made in 
(\ref{eq:generator}), we will, for each local field $K$ with the residue field
$k_K$, consider the 
collection ${\mathcal D}_K$ of additive characters $\psi:K\to \Cgordon^{\ast}$ 
satisfying
\begin{equation}\label{eq:D_K}
\psi(x)=\exp\left(\frac{2\pi i}{p}{\text{Tr}}_{k_K}(\bar x)\right)
\end{equation}
for $x\in {\mathcal O}_K$, where $p$ is the characteristic of $k_K$, 
$\bar x\in k_K$ is the reduction of $x$ modulo the uniformizer $\varpi_K$, and 
$\text{Tr}_{k_K}$ is the trace of $k_K$ over its prime subfield.
Any character from this collection can serve to produce an isomorphism between 
$K$ and $\hat K$. 
An example of a character from ${\mathcal D}_K$ is constructed in 
\cite[2.2]{Tategordon}. {It} is also naturally an exponential function.

\subsubsection{}One starts by formally adding exponential functions to the definable world. 
There are two kinds of exponentials one needs to add: the ones defined 
over the valued field, and the ones defined over the residue field.

For $Z$ in $\de_k$, the category $\rde_Z^{\text{exp}}$ consists of triples
$(Y\to Z,\xi,g)$, where $Y$ is in $\rde_Z$, and $\xi$, $g$ are morphisms in $\de_k$, $\xi:Y\to h[0,1,0]$, and $g:Y\to h[1,0,0]$. 
A morphism $(Y'\to Z,\xi',g') \to (Y\to Z,\xi,g)$ in  $\rde_Z^{\text{exp}}$
is a morphism $f:Y'\to Y$ in $\de_Z$ such that $\xi'=\xi\circ f$, and 
$g'=g\circ f$. 
The idea is that $e^{\xi}$ will be an exponential function on $Z$ 
over the residue field, and $e^g$ -- over the valued field.
We will soon define the Grothendieck ring $K_0(\rde_Z^{\text{exp}})$. 
The class $[Y\to Z,\xi, g]$ will be suggestively denoted by 
${\bf e}^{\xi}E(g)[Y\to Z]$.

Before we describe the relations that define the Grothendieck ring 
$K_0(\rde_Z^{\text{exp}})$, let us explain the intended specialization of 
constructible exponential functions to the $p$-adic fields.
As in Section~\ref{sec:back}, we will only consider 
${\mathcal L}_{\Zgordon}$-definable functions here. Recall that to interpret motivic
functions, we just needed to fix a field $K$ in ${\mathcal A}_{\Zgordon}$ or 
in ${\mathcal B}_{\Zgordon}$, and a 
uniformizer $\varpi_K$ of the valuation on $K$. 
To interpret exponential motivic {functions}, one needs in addition an element 
$\psi_K:K\to \Cgordon^{\ast}$ of the set 
${\mathcal D}_K$ of additive characters satisfying 
(\ref{eq:D_K}), as in Subsection~\ref{subsub:characters}.

Now suppose we have a triple $\varphi=(W,\xi,g)\in \rde_Z^{\exp}$, where $W$ 
is an ${\mathcal L}_{\Zgordon}$-definable subassignment equipped with an 
${\mathcal L}_{\Zgordon}$-definable morphism $\pi:W\to Z$, 
and $\xi, g$ -- ${\mathcal L}_{\Zgordon}$-definable morphisms 
from $W$ to $h[0,1,0]$ and 
$h[1,0,0]$, respectively.
For every $\psi_K$ in ${\mathcal D}_K$, we make a function 
$\varphi_{K,\psi_k}:Z_K\to \Cgordon$. Recall that the morphisms $\xi$ and $g$ give 
the functions $\xi_K:Z_K\to k_K$ and $g_K:Z_K\to K$ (all well defined when residue characteristic of $K$ is large enough). We define the function
$\varphi_{K, \psi_K}:Z_K \to \Cgordon$ by:
\begin{equation}\label{eq:interp}
z\mapsto 
\sum_{y\in \pi_K^{-1}(z)}\psi_K(g_K(y))
\exp\left(\frac{2\pi i}{p}\text{Tr}_{k_K}(\xi_K(y))\right).
\end{equation}

Now we are ready to define the Grothendieck ring $K_0(\rde_Z^{\text{exp}})$ that will play the same role as the ring $K_0(\rde_Z)$ played in the definition of constructible motivic functions in Section~\ref{sub:constr.f}.
The first relation is, as expected:
\begin{multline}\label{eq:rel1}
%\begin{aligned}
[(Y\cup Y')\to Z, \xi,g]+[(Y\cap Y')\to Z, \xi_{Y\cap Y'}, g_{Y\cap Y'}]\\
=[Y\to Z,\xi_Y,g_Y]+[Y'\to Z,\xi_{Y'},g_{Y'}].
%\end{aligned}
\end{multline}
for $Y,Y'\in \rde_Z$, and $\xi, g$ defined on $Y\cup Y'$.

The next relation is needed to take care of 
the restrictions of the exponential functions on the valued field to 
the residue field.
For a function $h:Y\to  k[[t]]$, denote by $\bar h$ its reduction 
$\mod (t)$, so that $\bar h:Y\to \A^1_k$. 
The second relation is:
\begin{equation}\label{eq:rel2}
[Y\to Z, \xi, g+h]=[Y\to Z, \xi+\bar h, g]
\end{equation}
for $h:Y\to h[1,0,0]$ a definable morphism with $\ord(h(y))\ge 0$ for all 
$y\in Y$.
Note that this condition becomes very natural in view of the interpretation (\ref{eq:interp}) and the condition (\ref{eq:D_K}) on the character $\psi_K$.

The third relation encompasses the fact that the integral of a character (of the residue field) over the field is zero. It postulates that
\begin{equation}\label{eq:rel3}
[Y[0,1,0]\to Z, \xi+p, g]=0
\end{equation}
when $p:Y[0,1,0]\to h[0,1,0]$ is the projection onto the second factor, and 
the morphisms $Y[0,1,0]\to Z$, $\xi$, and $g$ factor through the projection 
$Y[0,1,0]\to Y$. \footnote{Note that when $Y=Z=h_{\spec k}$ is a point, this statement literally amounts to the sum of the values of the character over the 
finite field being $0$. So, in general, this is {the} statement {that} the sum of the character over the {fibre} of $Y[0,1,0]$ over each point  $y\in Y$ equals $0$.}
  
Finally, the additive group of the Grothendieck ring $K_0(\rde_Z^{\exp})$ 
is defined as the group of formal linear combinations of equivalence classes of triples  $[Y\to Z, \xi, g]$ as above, modulo the subgroup generated by the relations (\ref{eq:rel1}), (\ref{eq:rel2}), and (\ref{eq:rel3}). 
It turns out that one can define multiplication on this set, so that the subgroup generated  by (\ref{eq:rel1}), (\ref{eq:rel2}), and (\ref{eq:rel3}) is an 
ideal \cite[Lem.~3.1.1]{CLF}, making  $K_0(\rde_Z^{\exp})$ into a ring.
This  ring is used instead of $K_0(\rde_Z)$ in the definition of 
constructible \emph{exponential} functions.

In \cite{CLF}, integration of constructible exponential functions is defined 
(along with the class of integrable functions), and that allows {one} to 
define {the} Fourier transform (satisfying all the expected properties).
The specialization principle holds 
for constructible {exponential} functions  as well, \cite[Th.~9.1.5]{CLF}. 
Namely, given a local field $K$,  
if we start with a constructible exponential function, integrate
it motivically, and then specialize the result to $K$
(using a fixed character $\psi\in {\mathcal D}_K$) 
according to the formula (\ref{eq:interp}), we would get the same result as if we had done the specialization (using the same character) first, 
and then integrated it with respect to the classical $p$-adic measure, 
when the residue characteristic $p$ is large enough.

\section{Appendix 1: the older theories}\label{geo}

Here we give very brief outlines of geometric motivic integration, 
and arithmetic motivic integration according to \cite{DL.arithm}, in order to point out the relationship of \cite{CL} with these older theories, 
and their relative features. In a sense, we are assuming some familiarity 
with geometric motivic integration, though the basic idea is sketched below.
There are excellent expositions \cite{blickle}, \cite{veys}.

\subsection{Arc spaces and geometric motivic measure}\label{sub:geo}
In the original theory of motivic integration,
the motivic measures live  on arc spaces of algebraic  varieties and 
take values in a certain  completion of the Grothendieck ring of the
category of all algebraic varieties over $k$.

Let $X$ be  a variety over $k$. 
The arcs are ``germs of order $n$ maps from the unit interval into $X$''.
Formally the space of arcs of order $n$ is defined as the scheme $\arcs_n(X)$
that represents  the functor
defined on the category of $k$-algebras by
$$R\mapsto \mathnormal{\mor_{k-{schemes}}}(\spec R[t]/t^{n+1}R[t],
X).$$
The {\bf space of formal arcs on $X$}, denoted by
$\arcs(X)$, is the inverse limit
$\lim\limits_{\longleftarrow}\arcs_n(X)$
in the category of $k$-schemes of the schemes $\arcs_n(X)$. 

The set of $k$-rational points of
$\arcs(X)$ can be identified with the set of points of $X$ over
$\ta{k}$, that is,
$$\mathnormal{\mor_{k-{schemes}}}(\spec\ta{k},X).$$
There are canonical morphisms 
$\pi_n : \arcs(X)\to \arcs_n(X)$ -- on the set of points, they correspond to
truncation of arcs.  In particular, when $n=0$,
we get the  the natural projection $\pi_X:\arcs(X)\to X$.

We use only the arc space of the $m$-dimensional affine space 
in these notes, so all that we need about arc spaces is 
essentially contained in the next example.

\begin{example} {(The arc space of the affine line $\arcs(\A^1)$.)}
By  definition, $\arcs_n(\A^1)$ represents
 the functor
\begin{align*}
R\to {\mor}(\spec R[t]/t^{n+1}R[t],\A^1)=\mor(k[x],R[t]/t^{n+1}R[t])\\
\cong
R[t]/t^{n+1}R[t]\cong R^{n+1}.
\end{align*}
Hence, $\arcs_n(\A^1)\cong\A^{n+1}$, and the natural projection
$\arcs_{n+1}(\A^1)\to\arcs_n(\A^1)$ corresponds to the map
$R[t]/t^{n+2}R[t]\to R[t]/t^{n+1}R[t]$ that takes
$P\in R[t]/t^{n+2}R[t]$ to $(P\pmod{t^{n+1}})$, which, in turn,
corresponds to the map $(T_0,\dots,T_{n+1})\mapsto(T_0,\dots,T_n)$
from $\A^{n+2}$ to $\A^{n+1}$.
We conclude that the inverse limit of the system $\arcs_n(\A^1)$ coincides with
the inverse limit of the spaces $\A^n$ with natural projections.
\end{example}

For simplicity, assume that the variety $X$ is smooth.  When $X$ is not 
smooth, the theory still works, but there is much more technical detail
(this constitutes the essence of \cite{DL}. See \cite{blickle} for 
an exposition).
A set ${\mathcal C}\subset \arcs(X)$ is  
called {\bf cylindrical} if it  
is of the form $\pi_n^{-1}(C)$ where $C$ is a constructible subset of
$\arcs_n(X)$. 
Let ${\mathcal C}$ be a cylinder 
with constructible base 
$C_n=\pi_n({\mathcal C})\subset \arcs_n(X)$. 
Then the motivic volume of ${\mathcal C}$ is by definition 
$\lef^{-n\dim(X)}[C_n]$, which is an element of $K_0(\var_k)[\lef^{-1}]$ 
(see Section~\ref{sub:gr.rings} for the definition of this ring). 

The \emph{geometric motivic measure} was initially defined as 
an additive function on an algebra
of subsets of the space $\arcs(X)$ that had a good approximation by cylindrical sets,  with values in 
a completion of $K_0(\var_k)[\lef^{-1}]$.

\subsection{Toward arithmetic motivic measure}

As we see from Section~\ref{sub:geo}, one could think of 
$\F_q[[t]]$-points of a 
variety as $\F_q$-points of its arc space. It is also possible to think of
$X(K)$ as $\arcs(X)(\F_q)$ for a characteristic zero local field 
$K$ with residue field $\F_q$. So the first idea would be to assign measures 
to cylindrical sets as above. However, there are two major problems with 
this approach. First, the the tools of the theory of geometric motivic integration that deal with singularities do not work when the residue field has finite characteristic, and restricting ourselves just to cylindrical subsets with a smooth base leaves us with far too few measurable sets. Of even greater 
importance is the issue that when the residue field is not algebraically 
closed, the action of the Galois group becomes {important}, and this Galois 
action varies with $p$. It turns out that Chow motives are ideally suited for keeping track of Galois action, and this is why arithmetic motivic measure takes values in a localization of $K_0(\mot_k)$ as opposed to $K_0(\var_k)$, which is 
sufficient when the field $k$ is algebraically closed.
For these reasons, the original theory of arithmetic motivic integration developed in \cite{DL.arithm} is quite different from geometric motivic integration, and is based more on logic than on algebraic geometry.

Fix the base field  $k$ of characteristic $0$ (for example, $k=\Qgordon$).
Let us first look at geometric motivic integration on $\A^m$ from the point 
of view of definable subassignments rather than arc spaces of varieties. 
The basic measurable sets for geometric motivic measure are stable cylinders in the arc space 
$\arcs(\A^m)$. 
Recall that a point in $\arcs(\A^m)$ can be thought of as an $m$-tuple of power series. 

Let ${\mathcal C}=\pi_n^{-1}(C_n)$ be a cylinder as in the previous subsection.
Suppose for simplicity that the set $C_n=\pi_n({\mathcal C})$ is defined just by one polynomial equation $f(\underline{x}^{(n)})=0$, where
$f(\underline {x}^{(n)})=f(x_1^{(n)},\dots, x_m^{(n)})$ is a polynomial with coefficients in $k$, and
$\underline{x}^{(n)}$ is an $m$-tuple of truncated power series, with  
each coordinate $x_i^{(n)}$, $i=1,\dots, m$ being a polynomial in $t$ of degree $n-1$.  
\begin{exercise}
The cylinder ${\mathcal C}$ is given by the following formula in the language of Denef-Pas:
$$
\phi_{\mathcal C}(\underline{x})=\exists \underline{y}\colon \ord(\underline{x}-\underline{y})\ge n\quad\wedge
\quad f(\underline{y})=0.
$$
Here $\underline{x}$ and $\underline{y}$ are $m$-tuples of variables  
ranging over
$k[[t]]$. 
\end{exercise} 

Thus, geometric motivic volume of the cylinder ${\mathcal C}$ is obtained in the following way:
the polynomial  $f(\underline{x})$,  (where $\underline{x}$ is an $m$-tuple of variables of the valued field 
sort) is replaced by the collection of ``truncated'' polynomials 
$f(\underline{x}^{(n)})$, 
where $\underline{x}^{(n)}$ is now an $m$-tuple of polynomials in $t$ 
of degree $n-1$ with coefficients in the residue field sort (\ie, in $k$).
In  the next step, the condition   $f(\underline{x}^{(n)})=0$ is replaced by the collection of equations stating that
all the coefficients of the resulting polynomial in $t$ equal zero, which defines a constructible set over the 
residue
field. Finally, the motivic volume of ${\mathcal C}$ is the class of the constructible set obtained above
multiplied by $\lef^{-nm}$, where $n$ is the 
level of truncation.

If we state the basic idea of geometric motivic integration in this form, it becomes natural to 
define the corresponding procedure for more general formulas in Denef-Pas 
language, than just the formulas defining  cylinders. 
The following two key steps allowed the above construction to work:
first, we were able to replace the formula that had a quantifier over the valued field 
with a formula without quantifiers over the valued field; and then 
the value of the motivic volume was obtained  from a ring formula with 
variables 
in the residue field (the formula defining the constructible set 
$\pi_n({\mathcal C})$).
In general, both of these steps rest on the process of quantifier elimination -- that is, replacing 
a formula that has quantifiers with an equivalent formula with no quantifiers (see \cite{Tom.intro}
for a discussion of quantifier elimination in this context).

\subsection{Quantifier elimination}

The following theorem which is a version of the theorem of Pas allows to eliminate all
quantifiers in DVR-formulas except the ones over the residue field.

\begin{theorem} 
Suppose that $R$ is a ring of characteristic zero.
Then for any DVR-formula $\phi$ over $R$ there exists a DVR-formula $\psi$ over $R$  which contains no quantifiers running over the valuation ring and 
no quantifiers running over $\Zgordon$, such that:
\begin{enumerate}
\item $\theta\leftrightarrow\phi$ holds in $K[[t]]$ for all fields $K$ containing $R$,
\item $\theta\leftrightarrow\phi$ holds in $\Zgordon_p$ for all $p\gg0$ when $R=\Zgordon$.
\footnote{Even if the original formula $\phi$ had no {quantifiers} over the residue field, the formula 
$\psi$ might have them.}
\footnote{Elimination of quantifiers over the value sort is 
due to Presburger.
It is  because we want this quantifier elimination  result to hold, multiplication is not
permitted for variables of the value sort (it is the famous {theorem of G\"{o}del} that $\N$ with the 
standard operations 
does not admit quantifier elimination).
}
\end{enumerate}
\end{theorem}

\begin{theorem}(Ax)\cite[\S~8.2]{FJ}
Algebraically closed fields admit elimination of quantifiers in the language of rings.
\end{theorem}

In particular, this theorem implies the theorem, due to Chevalley, 
stating that an image of a constructible set under a projection morphism is 
a constructible set.

The situation is different for non-algebraically closed fields; in particular,
 the quantifiers over the residue 
field of a local field cannot be eliminated in general.\footnote{A theorem due to A. Macintyre \cite{Mac} states that there would be complete quantifier elimination
if we added to the language, for each $d$, the predicate ``$x$ is the $d$-th power in the field''. 
So in some sense all quantifiers except in the formulas '$\exists y: y^d=x$' can be eliminated.
This fact is reflected in the theory of Galois stratifications, which is the main tool in the construction of the map that takes the motivic volume of a 
definable set into the Grothendieck ring of the {category} of Chow motives, 
discussed in Section~\ref{sub:comp}.} 

Ax's theorem is, in some sense, the reason why geometric motivic measure
is so much easier to construct than arithmetic motivic measure. In the easiest case of the stable 
cylinder, for example, once we have the ring formula over the residue field, quantifier elimination 
produces a quantifier-free formula over the residue field, that is, a constructible set.

\subsection{The original construction of arithmetic motivic measure} 
The original construction of {arithmetic} motivic {measure}
\cite{DL.arithm} {follows} these steps.
\begin{enumerate}
\item[0.] We start with a DVR-formula $\phi$ or, equivalently, with a definable subassignment of the functor
$h_{\arcs({\A_m})}$. When interpreted over a $p$-adic field, the formula $\phi$ gives a measurable 
set (in the classical sense).   
\item[1.] For every positive integer $n$, the definable subassigment $h$ 
defined by $\phi$ can be truncated at level $n$. By Pas's theorem on elimination of quantifiers,
the truncated subassignment $h_n$ is definable by a ring formula $\psi_n$ over the residue field 
(note that the number of variables of $\psi_n$ depends 
on $n$).\footnote{
By  
analogy with the corresponding notion in the construction of  Lebesgue measure, one can say that 
the formulas $\psi_n$ define the ``outer'' approximations to the set defined by the formula 
$\phi$ (see \cite{Tom.intro}).}
\item[2.] We consider the class of the formula $\psi_n$ in $K_0(\pff_k)$. At this step, essentially, 
the formulas that should have the same motivic volume are getting identified.
\item[3.] There is a map from $K_0(\pff_k)$ to the Grothendieck 
ring of the category of Chow motives. One takes the virtual Chow motive $M_n$ 
associated with $[\psi_n]$. 
\item[4.] In the original construction of \cite{DL.arithm}, the ring of virtual Chow motives is completed
in a similar way to the completion of the ring of values of the geometric motivic measure.
Finally, the Chow motive associated with the definable subassignment $h$ is the inverse limit of $M_n\lef^{-nd}$,
where $d$ is the dimension of $h$.
\end{enumerate}

\begin{remark}
As we have seen from Section~\ref{sub:geo}, the construction of geometric motivic measure follows the same steps, 
with the following simplifications: In Step 2, we need to consider 
{the equivalence} relation on formulas that
comes not from comparing them on pseudofinite fields, but from comparing them on algebraically closed fields.
Instead of the complicated Step 3, one can apply quantifier elimination over algebraically closed fields to it
to obtain a class of a constructible set, that is, an element of $K_0(\text{Var}_k)$.
\end{remark}

\subsection{Measurable sets in different theories}

It is worth pointing out that almost every variant of motivic 
integration has a 
slightly different algebra of measurable sets. 
In the very first papers on motivic integration, {\eg} \cite{DL.McKay}, 
the measurable sets 
were the semi-algebraic sets, and then later $k[t]$-semi-algebraic sets.
In geometric motivic integration that developed later, 
the basic measurable sets are 
stable cylinders. In \cite[Appendix]{DL.McKay}, a measure theory and a 
$\sigma$-algebra of measurable sets that includes stable cylinders 
is worked out. 
Looijenga \cite{looijenga} describes a slightly different version of 
geometric motivic integration as well. 
Here we make a few remarks about the relationships between the algebras of measurable sets in all these articles.

\subsubsection{Cylinders {\it vs.} semi-algebraic sets}
The algebra of sets definable in Denef-Pas language with coefficients in 
$k[t]$ specializes to the algebra of the $k[t]$-semi-algebraic sets.
It follows from J. Pas's theorem on quantifier elimination that
if the set $A$ is semi-algebraic, then $\pi_n(A)$ is a constructible 
subset of ${\mathcal L}_n(X)$.
This statement ultimately implies that the algebra of semi-algebraic subsets is contained in the
algebra of measurable sets of \cite[Prop.~1.7 (2)]{DL.McKay}.
The advantage of working with measurable sets that are well approximated by cylinders is that this algebra is more geometric, and bigger than the algebra of $k[t]$-semi-algebraic sets. The main disadvantage is that one needs to complete the ring $K_0(\var_k)$ in order to define the measure on this algebra.
On the other hand, the algebra of $k[t]$-semi-algebraic sets possesses two 
advantages:
first, it is this algebra that we can get by specializing the theory of Cluckers and Loeser to algebraically closed fields, and so it follows that it is 
not necessary to complete
the ring $K_0(\var_k)$ in order to define the restriction of the 
motivic measure to this algebra -- inverting $\lef$ and $1-\lef^{-n}$, $n>0$ 
is sufficient. 
Second, it is this algebra that (at the moment) appears in all the 
generalizations 
of the motivic measure theory (\ie, in motivic integration on formal schemes \cite{Sebag}).

\subsubsection{Denef and Loeser {\it vs.} Looijenga}
Denef and Loeser work directly with subsets of the underlying topological space of the 
$k$-scheme ${\mathcal L}(X)$, whereas Looijenga considers 
subsets of the space of sections of its {structure} 
morphism (as a scheme over $\spec k[[t]]$) which is in bijection with the 
set of \emph{closed} points
of ${\mathcal L}(X)$. 
Thus, the algebra of measurable subsets constructed in \cite{looijenga}
is the restriction of the algebra of measurable sets of \cite{DL} 
to the set of closed points
of ${\mathcal L}(X)$.
The advantage of the approach in \cite{looijenga} 
is that it makes no difference between the schemes $X$ originally defined over $k$ {\it vs.} the schemes
$X$ defined over $R=k[[t]]$, which is sometimes very useful in applications. 
Indeed, if $X$ is defined over $k$, it can be base changed to $k[[t]]$: 
set ${\mathcal X}:=X\times_{\spec k}\spec k[[t]]$,
and then the set $\mx_{\infty}$ of \cite{looijenga}, which is in bijection with the set of closed points of $\arcs(X)$, is defined as the set of sections of the structure 
morphism of the scheme $\mx$.  

\section{Appendix 2: an example}

Here we do in detail a calculation of the motivic volume of a set that is 
relevant to character values of depth zero representations of 
$G=\text{SL}(2,K)$ for a local field $K$. 
Complete character tables for depth zero representations of $\text{SL}(2,K)$ 
restricted to the set of topologically unipotent elements appear in 
\cite{CG}, and we refer to that article for the detailed explanation {as to} why these 
sets appear.

Roughly, the calculation goes as follows. Depth zero representations are
obtained from representations of finite groups by inflation to 
a maximal compact subgroup followed by compact induction.
There is a well-known Frobenius formula for the character of an induced 
representation, and it applies in this situation as well. Let $H$ be a compact 
set of topologically unipotent elements of $G$, let $f_H$ be the characteristic function of this set. Let $\pi$ be a depth zero representation that is induced from the maximal compact subgroup $G_x=\text{SL}(2, \rigordon_K)$, and let 
$\Theta_{\pi}$ be its Harish-Chandra character.
Then (see \cite{CG})
$$
\Theta_{\pi}(f_H)=\mu(G_x)\int_{G/G_x}\int_H\chi_{x,0}(g^{-1}hg)\,dh\,dg,
$$
where $\chi_{x,0}$ is the inflation to $G_x$ of the character of the 
representation of $\text{SL}(2,\F_q)$ that $\pi$ restricts to.

It is, therefore, natural that the volume of the set of elements $g\in G$ 
such that the element $g^{-1}hg$ is in $G_x$ and projects under 
reduction $\mod \varpi$ 
to a given unipotent conjugacy class of $\text{SL}(2,\F_q)$ is the key to 
the value of the character at $h$. 

The following calculation appears when we take $h$ of the form
$h=\left[\begin{smallmatrix} 0 & \varpi^n u\\ 
\epsilon\varpi^{n} u& 0\end{smallmatrix}\right]$, and the unit $u$ is a square.

Even though this calculation is included as an explicit example of a computation of a motivic volume, we are not using the technique of inductive application
of cell decomposition (since there are four valued field variables, it would have been too tedious a process). Instead, we do the calculation {\it ad hoc}, using the  older approach through the outer motivic measures that is sketched in the previous appendix. 
However, our motivic volume depends on a residue-field parameter, so we are using the language and the results of \cite{CL} as well. 

Let us first introduce an abbreviation for the 
``reduction $\mod \varpi$'' map: let
$$
\bar x=\begin{cases}
\ac(x), & \ord(x)=0\\
0,& \ord(x)>0.
\end{cases}
$$

\begin{example}
Let us consider the family of formulas depending on a parameter 
$\eta$ that ranges over the set of non-squares in $\F_q$ (note that this
is a definable set):
\begin{equation}
\phi_{\eta}(a,b,c,d)\quad=\text{`}ad-bc=1 \quad\wedge\exists\, \xi
(\bar b^2-\bar d^2\eta=\xi^2)\text{'}.
\end{equation}
We claim that the motivic volume of the set defined by $\phi_{\eta}$ is independent of $\eta$ and equals $\frac12\lef(\lef-1)(\lef+1)$.

\begin{proof}
Consider the formula 
%\begin{multline}
\begin{equation}
%\begin{aligned}
\Phi(a,b,c,d,\eta)
=\text{`}
ad-bc=1\ \wedge
\exists\, \xi \neq 0(\bar d^2-\bar b^2\eta=\xi^2)
\ \wedge
\nexists\, \beta (\eta=\beta^2) 
\text{'}.
%\end{aligned}
\end{equation}
%\end{multline}
In this formula, four variables $a,b,c,d$  range over the valued
field, the variable $\eta$ ranges over the residue field, and 
all the quantifiers range over the residue field.
It defines a subassignment of $h[4,1,0]$, which corresponds to the
disjoint union of the subassignments of $h[4,0,0]$ defined by the formulas
$\phi_{\eta}$ over all non-squares $\eta$. 

%\begin{mysteps}
\subsubsection{Step 1. Reduction to the residue field}
The formula $\Phi$ can be broken up into two parts according to whether 
$b$ is a unit: $\Phi=(\Phi\wedge (\ord(b)=0))\vee (\Phi\wedge (\ord(b)>0))$.

We start by showing that the subassignment defined by $\Phi\wedge
(\ord(b)=0)$ is
stable at level  $0$, \ie, that it is essentially ``inflated'' from
the finite field. In order to do this, we need to introduce an
abstract variety $V$ that plays the role of the ``projection'' of this
formula to the residue field.

Let $k$ be an arbitrary field of characteristic zero 
(the theory of arithmetic motivic integration tells us that we should
think of  $k$ as a pseudofinite field).
Consider  the   subvariety  $V$ of $\A^4$ over $k$ cut out 
by the equation $x_1^2-x_2^2x_3=x_4^2$. It has no singularities
outside the hyperplane $x_2=0$ (note that this statement is true in any
characteristic greater than $2$). 
Recall the notation $[\phi]$ for the motivic volume of a formula
$\phi$ in the language of rings. Let
\begin{equation}\label{eqn:M1}
{\mathbb M}_1:=[\text{`}(x_1^2-x_2^2x_3=x_4^2)\wedge (x_2\neq 0)\wedge (x_3\neq 0)\wedge (x_4\neq
0)\wedge 
(\nexists \beta ~ (x_3=\beta^2))\text{'}].
\end{equation}
Consider the formula
%$$
\begin{multline}
\Phi_1(b,d,\eta,\xi):=\\
\text{`}
(\bar d^2-\bar b^2\eta=\xi^2)\wedge (\xi\neq 0)\wedge
(\ord(b)=0)\wedge (\eta\neq 0)\wedge(\nexists\,\beta~(\eta=\beta^2)
)\text{'}
\end{multline}
%$$
Set $x_1=\bar d$, $x_2=\bar b$, $x_3=\eta$, $x_4=\xi$. This
``reduction'' takes the formula $\Phi_1$ exactly to the ring formula  
that appears in the right-hand side of (\ref{eqn:M1}). Since it is
mutually exclusive with the formula $x_2=0$ that defines a set
containing the singular locus of $V$, 
the subassignment of $h[2,2,0]$ defined by the formula $\Phi_1$    
is stable at level $0$, and its motivic volume equals ${\mathbb M}_1$. 

Let 
$\Phi_2(b,d,\eta)$ be the formula 
$$
\text{`}\exists\xi~ 
(\bar d^2-\bar b^2\eta=\xi^2)\wedge (\xi\neq 0)
\wedge (\ord(b)=0)\wedge(\eta\neq 0)\wedge(\nexists \beta (\eta=\beta^2))
\text{'}.
$$
 
The formula $\Phi_1$ is a double cover of $\Phi_2$, so 
$\mu(\Phi_2)=\frac12\mu(\Phi_1)=\frac12{\mathbb M}_1$.

Finally,  consider 
the projection  $(d,b,c,a,\eta)\to(d,b,\eta)$.
The subassignment defined by $\Phi\wedge (\ord(b)=0)$ 
projects to the subassignment
defined by $\Phi_2$, and the volume
of each fibre of this projection is $\lef$: indeed, given that $b$ is
a unit, for every value of $a$, there is unique $c$ such that $ad-bc=1$. 
Hence, we have  $\mu(\Phi\wedge (\ord(b)=0))=\frac 12\lef{\mathbb M}_1$. 
%\end{mysteps}

%\begin{mysteps}
\subsubsection{Step 2. Independence of the parameter $\eta$}
This step consists in the observation that for each 
$\eta_1, \eta_2\in K^{\ast}\setminus {K^{\ast}}^2$, there is a
definable bijection between the triples $(x_1,x_2,x_4)$ and
$(x_1',x_2',x_4')$ 
such that 
$x_1^2-x_2^2\eta_1=x_4^2$ and ${x_1'}^2-{x_2'}^2\eta_2={x_4'}^2$: the {bijection}
is defined by the formula 
%\begin{equation*}
\begin{multline*}
\Psi(x_1,x_2,x_4,x_1',x_2',x_4',\eta_1,\eta_2)\\
=\text{`}\exists\,
y (\eta_1=\eta_2y^2) \wedge (x_1'=x_1) \wedge (x_2'=x_2y) 
\wedge (x_4'=x_4)\text{'}
 \wedge \nexists\, \beta (\eta_1=\beta^2).
\end{multline*}
%\end{equation*}

The Corollary \cite[14.2.2]{CL} together with
Remark \cite[14.2.3]{CL} implies 
that if the motivic volume is constant on the fibres, then the total
volume is the volume of the fibre times the class of the base.
It follows that for each $\eta\in\F_q\setminus{\F_q}^2$, 
we have 
$$
\mu(\phi_{\eta}\wedge (\ord(b)=
0))=\frac2{\lef-1}\mu(\Phi\wedge(\ord(b)=0)).
$$
%\end{mysteps}

%\begin{mysteps}
\subsubsection{Step 3. A residue-field calculation: finding ${\mathbb M}_1$}
Note that it is in this step that we see the conic promised in the
introduction.

We start by considering abstract varieties again. Recall the variety
$V$ defined in Step 1. Let us denote the coordinates on $\A^3$ by
$(t,s,e)$, and consider the subvariety $V_2$ of $\A^3$ defined by 
the equation $t^2-s^2=e$.
Consider the birational map 
$(x_1,x_2,x_3,x_4)\mapsto(x_2,\frac{x_1}{x_2},x_3,\frac{x_4}{x_2})$ from the
variety $V$  to the variety $V_2\times\A^1$. It is an
isomorphism  between the open sets $x_2x_3\neq 0$ in $V$ and
$e\neq 0$ in $V_2\times \A^1$.
Then the class ${\mathbb M}_1$ equals $[\text{`}\nexists\beta\neq 0
  :t^2-s^2=\beta^2\text{'}](\lef-1)$. It remains to calculate the class  
$[\nexists\beta\neq 0 :t^2-s^2=\beta^2]$.

The class $\lef^2$ of the $(t,s)$-plane breaks up into the sum of the three classes:  
$$
\lef^2=
[\nexists\beta (t^2-s^2=\beta^2)]+
[\exists \beta (t^2-s^2=\beta^2)\wedge \beta\neq 0]+[t^2-s^2=0].
$$
It is easy to see that $[t^2-s^2=0]=2(\lef-1)+1$.
We also have
$[\exists \beta (t^2-s^2=\eta^2)\wedge \beta\neq 0]=\frac12(\lef-1)[x^2-y^2=1]
=\frac12(\lef-1)(\lef-1)$. 
Therefore, $[\nexists\beta (t^2-s^2=\eta^2)]=\frac12\lef^2-\lef+\frac12$.
Hence, ${\mathbb M}_1=\frac12(\lef-1)^3$.

Finally, we have:
\begin{equation}\label{eq: the hard volume}
\begin{aligned}
\mu(\phi_{\eta}\wedge \ord(b)=0))&=\frac2{\lef-1}\mu(\Phi\wedge (\ord(b)=0))\\
&=\frac2{\lef-1}\frac14(\lef-1)^3\lef=\frac12\lef(\lef-1)^2.
\end{aligned}
\end{equation}
%\end{mysteps}

%\begin{mysteps}
\subsubsection{Step 4. Completing the proof}  
It is easy to calculate the motivic volume of the remaining part 
 $\phi_{\eta}\wedge (\ord(b)>0)$.
If $\ord(b)>0$, then the formula $\psi(\bar d,\bar b,\eta)$ becomes 
'$\exists \beta\neq 0 (\bar d^2=\beta^2)$'.
Clearly, its motivic volume  is $(\lef-1)$.
It remains to notice that if $\ord(b)>0$, 
then the variable $c$ contributes the factor $\lef$, 
both $(a,d)$  have to be units, and once $d$ is chosen, $a$ is
 determined uniquely by the determinant condition. 
Altogether, we get  
$\mu(\phi_{\eta}\wedge (\ord(b)>0))=\lef(\lef-1)$.
Finally, we get:
$$
\begin{aligned}
\mu(W_{U_{\epsgordon},n/2}^{(0)}(h))&=\mu(\phi_{\eta})=
\frac12\lef(\lef-1)^2+\lef(\lef-1)\\
&=\frac12\lef(\lef-1)(\lef+1),
\end{aligned}
$$
which completes the proof.
%\end{mysteps}
\end{proof}
\end{example}

\end{document}